\documentclass[preprint,12pt]{elsarticle}

\usepackage{graphicx}
\usepackage{amsmath,amssymb}
\usepackage{amssymb}
\usepackage{amssymb}
\usepackage[version=4]{mhchem}
\usepackage{geometry}
\usepackage{enumitem} 
\usepackage{lineno}
\usepackage[bookmarks,bookmarksnumbered]{hyperref}
\hypersetup{colorlinks = true,linkcolor = blue,anchorcolor =red,citecolor = blue,filecolor = red,urlcolor = red}

\journal{}

\begin{document}

\begin{frontmatter}


\title{The Quasi-Steady-State Approximations
revisited: Timescales, small parameters, singularities, and normal forms
in enzyme kinetics}



\author[label1]{Justin Eilertsen}

\author[label1,label2]{Santiago Schnell\fnref{cor1}}
\ead{schnells@umich.edu}

\address[label1]{Department of Molecular \& Integrative Physiology, University of Michigan 
Medical School, Ann Arbor, MI 48109, USA}
\address[label2]{Department of Computational Medicine \& Bioinformatics, University of 
Michigan Medical School, Ann Arbor, MI 48109, USA}

\fntext[cor1]{Corresponding author}

\begin{abstract}
In this work, we revisit the scaling analysis and commonly accepted conditions for 
the validity of the standard, reverse and total quasi-steady-state approximations through 
the lens of dimensional Tikhonov--Fenichel parameters and their respective 
critical manifolds. By combining Tikhonov--Fenichel parameters with scaling 
analysis and energy methods, we derive improved upper bounds on the approximation error for the 
standard, reverse and total quasi-steady-state approximations. Furthermore, previous 
analyses suggest that the reverse quasi-steady-state approximation is 
only valid when initial enzyme concentrations greatly exceed initial substrate
concentrations. However, our results indicate that this approximation can be 
valid when initial enzyme and substrate concentrations are of equal magnitude. Using energy methods, we find that the 
condition for the validity of the reverse quasi-steady-state approximation is far 
less restrictive than was previously assumed, and we derive a new ``small'' 
parameter that determines the validity of this approximation. In doing so, we extend the established domain of validity for the reverse quasi-steady-state approximation. 
Consequently, this opens up the possibility of utilizing the reverse 
quasi-steady-state approximation to model enzyme catalyzed reactions and 
estimate kinetic parameters in enzymatic assays at much lower enzyme to substrate 
ratios than was previously thought. Moreover, we show 
for the first time that the critical manifold of the reverse quasi-steady-state
approximation contains a singular point where normal hyperbolicity is lost. 
Associated with this singularity is a transcritical bifurcation, and the 
corresponding normal form of this bifurcation is recovered through scaling 
analysis.
\end{abstract}


\end{frontmatter}


\section{Introduction}
Perhaps the most well-known reaction in biochemistry is the Michaelis-Menten (MM) reaction mechanism, (\ref{eq:mm1}), which describes the catalytic conversation of a substrate, $S$, into a product, $P$. Catalysis is achieved through means of an enzyme, $E$, that reversibly binds with the substrate, forming an intermediate complex, $C$. In turn, $C$ irreversibly disassociates into enzyme and product molecules:
\begin{equation}\label{eq:mm1}
    \ce{$S$ + $E$ <=>[$k_1$][$k_{-1}$] $C$ ->[$k_2$] $E$ + $P$ }.
\end{equation}

For theorists and applied mathematicians, whose aim is to mathematically describe the dynamics of chemical networks and metabolic pathways \cite{2004-Crampin-PBMB2,SAURO20045}, the MM reaction mechanism serves as a building block to these more complex systems. In general, there are two ways to mathematically model enzyme catalyzed reactions. At low concentrations of chemical species stochastic models are generally favorable since they describe the seemingly random collisions of reactant molecules in intracellular environments~\cite{TURNER2004165,Grima}. In contrast, when concentrations are high and the chemical species are well-mixed, the MM mechanism can be appropriately modeled as a set of nonlinear ordinary differential equations. Although one can argue that deterministic models are in some sense more manageable than stochastic models, nonlinear deterministic models rarely admit closed form solutions, and therefore model reduction techniques must be employed in order to simplify the model equations so that approximate solutions can obtained and analyzed. Typically, model reduction is synonymous with approximating the model dynamics on an invariant manifold, and the advent of powerful computer algebra systems aids in the systematic reduction of high-dimensional or even infinite-dimensional (i.e., partial differential equations and delay differential equations) dynamical systems. Not surprisingly, there is a large body of literature that clearly illustrates how model reduction can uncover, quantify, and explain the various nonlinear phenomena arising not only in chemistry, but also in biology \cite{BERTRAM2017105,Rubin, BOOTH201454} and physics \cite{Spiegel1983, ROBERTS1997,roberts1989,Crawford}. 

Historically, the most widely-utilized reduction technique in deterministic enzyme kinetics has been the singular perturbation method \cite{Omalley}, in which a reduced model is constructed by approximating the flow of the model equations on a \textit{slow invariant manifold} (SIM). Since the dimension of the SIM is less than the dimension of the phase-space, approximation of the dynamics on the SIM permits a \text{reduction} in the dimension of the problem. The singular perturbation method exploits the presence of disparate fast and slow timescales; the work of Tikhonov \cite{tikhonov1952} and Gradshtein \cite{Gradshtein1953} provides the necessary rigorous foundation for the construction of a reduced model. Briefly, when fast and slow timescales are present, the differential equations that model 
the MM reaction with a timescale separation can be expressed in the form
\begin{subequations}
\begin{align}
\dot{x} &= f(x,y;\varepsilon),\label{eq:1}\\
\varepsilon \dot{y} &= g(x,y;\varepsilon),\label{eq:2}
\end{align}
\end{subequations}
where $0 < \varepsilon \ll 1$. Setting $\varepsilon =0$ yields a
differential-algebraic-equation
\begin{subequations}
\begin{align}
\dot{x} &= f(x,y;0)\label{eq:DAE1},\\
0 &= g(x,y;0)\label{eq:DAE2},
\end{align}
\end{subequations}
and, according to Tikhonov's theorem, equations (\ref{eq:DAE1})--(\ref{eq:DAE2}) provide a very good approximation to the dynamics of (\ref{eq:1})--(\ref{eq:2}) when $\varepsilon$ is sufficiently small. The reduced system (\ref{eq:DAE1})--(\ref{eq:DAE2}) provides a simpler, and often times more tractable, \textit{reduced} 
mathematical model that is commonly referred to as a \textit{quasi-steady-state approximation} (QSSA). 

The hope is that the condition that supports the validity of (\ref{eq:DAE1})--(\ref{eq:DAE2}) (i.e., $\varepsilon \ll 1$) is easy to implement in enzymatic assays, so that precise and accurate measurements of the kinetic parameters pertinent to the reaction can be made by fitting the QSSA model to experimental time course data \cite{2016-Stroberg-BPC,KIM}. Due to the necessity of disparate timescales, slow manifold reduction may not be possible in every physical scenario. Therefore, the challenge for theorists is not only to derive a reduced model that has suitable utility, but also to determine the unique physical and chemical conditions that \textit{permit} the validity of the associated reduction. Thus, an important task of the theoretician is determine the exact conditions for which the reduced model is valid~\cite{roberts2018backwards}. Mathematically, this translates to determining ``$\varepsilon$," a (typically) dimensionless parameter that
may not be unique. The non-uniqueness of ``$\varepsilon$" adds complication, since  some ``epsilons'' are better than others.  The best-known example of the ``non-uniqueness dilemma" resides the history of the 
derivation of the Michaelis--Menten (MM) equation
\begin{equation} \label{eq:MM}
v = \frac{V s}{K_M+s},
\end{equation}
which is obtained by applying the QSSA to the MM reaction mechanism~(\ref{eq:mm1}).
In (\ref{eq:MM}), $v$ is the velocity of product formation in the reaction, 
$V$ is the limiting rate of the reaction, $K_M$ is the Michaelis constant, and 
$s$ is the free substrate concentration for the reaction. Alternatively, the MM equation is often referred to as the \textit{standard} quasi-steady-state approximation (sQSSA). In 
1967, Heineken et al.~\cite{HEINEKEN196795} formally applied, for the first time, 
the standard QSSA to the nonlinear differential equations governing the MM reaction 
mechanism~(\ref{eq:mm1}) via singular perturbation analysis. Based on the findings 
of Laidler~\cite{Laidler:1955:TTK}, a pioneer in chemical kinetics and authority 
on the physical chemistry of enzymes, Heineken et al.~\cite{HEINEKEN196795} 
introduced a consistent ``$\varepsilon$'' for the MM reaction mechanism through 
scaling analysis. A more rigorous analysis of the sQSSA was introduced by
Reich and Sel'kov~\cite{REICH1974} and Schauer and Heinrich \cite{SCHAUER1979A}, 
from which other ``epsilons'' where determined by proposing conditions that minimized the errors in the implementation of the sQSSA. In 1988, Lee A. Segel~\cite{Segel1988} derived 
the widely accepted criterion for the validity of the sQSSA and the derivation of 
the MM equation~(\ref{eq:MM}) by estimating the slow and fast timescales of the 
reaction. Segel~\cite{Segel1988} illustrated that prior physico-chemical knowledge 
about the reaction dynamic is instrumental in deriving the fast and slow timescales 
and uncovering the most general criteria for the validity of the sQSSA. As a direct 
result of Segel's scaling method, the conditions for the validity of the sQSSA were 
derived for suicide substrates~\cite{BURKE199081,Burke199397}, 
alternative substrates~\cite{2000-Schnell-JMC}, fully competitive enzyme 
reactions~\cite{2000-Schnell-BMB1}, zymogen activation~\cite{BPC}, and coupled 
enzyme assays~\cite{MBS,JTB}. Segel's scaling 
approach was also applied to the analysis of the MM reaction mechanism~(\ref{eq:mm1}) to extend the 
validity of the QSSA in different regions of the initial enzyme and 
substrate concentration parameter space via the reverse QSSA
(rQSSA)~\cite{Segel1989,SCHNELL2000}, and the total QSSA
(tQSSA)~\cite{BORGHANS1996,SCHNELL2002,Tzafriri2003}.

Despite the power of Segel's scaling and simplification analysis for the MM 
reaction mechanism~(\ref{eq:mm1}), there is still a fundamental challenge with 
its implementation via the rQSSA: there has never been a small parameter (i.e., a 
specific ``$\varepsilon$''), analogous to the one obtained by Segel~\cite{Segel1988} 
for the sQSSA, that is as effective at determining when the rQSSA is valid. Unfortunately, estimating the fast timescale associated with the rQSSA is difficult, and there have been fundamental disagreements in the reported estimates~\cite{SCHNELL2000}. Thus, 
timescale estimation continues to be the ``Achilles' heel'' of the rQSSA. This 
raises the obvious question of whether or not timescale estimation is truly the 
best approach towards resolving this problem.

Recently, Walcher and his collaborators~\cite{GOEKE1,Goeke2014,Noethen2011} demonstrated that identifying a
Tikhonov--Fenichel parameter (TFP) is an effective way to determine a priori 
conditions that suggest the validity of the QSSA. Essentially, a TFP is a dimensional parameter -- such as a rate constant or initial concentration of 
a species -- that, when identically zero, ensures the existence of a manifold of equilibrium points. Such manifolds are central to \textit{geometric singular perturbation theory} (GSPT) and, as a result of Fenichel \cite{GSPT0}, it is well-understood that when certain conditions hold, the existence of a critical manifold of equilibria ensures the existence of a SIM once the TFP is small but non-zero. In this sense the \textit{origin} of the SIM can be linked to the vanishing of a specific dimensional parameter, and a sufficiently small TFP ensures the existence of a SIM and the corresponding validity of a QSSA. However, the identification of a TFP does not diminish the importance of the asymptotic small parameter $\varepsilon$, 
since it is essential to define what physically constitutes ``small'' when a 
parameter is non-zero. In other words, we must still answer the question: 
\textit{how} small should Tikhonov--Fenichel parameters be in comparison 
to other dimensional parameters in order to yield an accurate reduced model? 

In this paper, our primary objective is to determine a specific small parameter 
that determines the validity of the rQSSA to the MM reaction mechanism~(\ref{eq:mm1}), but also to convey subject matter that can be quite technical in a language capable over reaching a wider audience that is not limited to applied mathematicians and physical chemists.
In Section \ref{sec:revisit}, we review the conditions for the validity of the various
quasi-steady-state approximations that are commonly employed to approximate the 
long-time dynamics of the MM reaction mechanism~(\ref{eq:mm1}): namely, the 
sQSSA, the rQSSA, and tQSSA. In Section~\ref{sec:GSPT}, we assess the 
validity of specific QSS reductions by employing geometric singular perturbation 
theory, and illustrate how assumptions about the validity of the sQSSA based on 
Segel's timescale separation can lead to erroneous conclusions. In 
Section~\ref{sec:nonscaled}, we introduce two methods that do not rely on scaling 
or timescale separation: Tikhonov--Fenichel parameters and energy methods, both of which can be employed to determine the validity of 
the QSSA. In Section~\ref{sec:rQSSA}, we analyze the validity of the rQSSA using 
the methods discussed in Section~\ref{sec:nonscaled}, and in Section~\ref{sec:timescales}, we discuss timescale separation and the \textit{hierarchy} of small parameters that support the justification of the sQSSA, rQSSA, and the tQSSA. Finally, in our discussion
(Section~\ref{sec:discussion}), we summarize our results and critique some of the 
conclusions drawn about the validity of the sQSSA and the rQSSA in the 
previous analyses of Segel and Slemrod~\cite{Segel1989}, and Schnell and 
Maini~\cite{SCHNELL2000}, respectively. We also discuss the impact our results 
will have on experimental assays, and how the methods we utilize can be employed 
to analyze more complicated reactions and experimental assays.   

\section{Applying the Quasi-Steady-State Approximations to the Michaelis--Menten 
reaction mechanism: Scaling and simplification approaches} \label{sec:revisit}
In this section we review the application of the different versions of the QSSA
for the MM reaction mechanism~(\ref{eq:mm1}), and the mathematical justification for 
the application of each approximation. We also present the timescales 
and criterion for the validity of the QSSA originally derived by Segel~\cite{Segel1988} 
and Segel and Slemrod~\cite{Segel1989}.

\subsection{Heuristic estimation of fast and slow timescales using the 
\textit{standard} Quasi-Steady-State Approximation}
The mathematical model that describes (deterministically) the reaction 
mechanism~(\ref{eq:mm1}) is a set of nonlinear ordinary differential 
equations,
\begin{subequations}
\begin{align}
\dot{s} &= -k_1(e_0-c)s+k_{-1}c\label{eq:1a},\\
\dot{c} &=  k_1(e_0-c)s -(k_{-1}+k_2)c\label{eq:1b},\\
\dot{p} &= k_2c\label{eq:1c}.
\end{align}
\end{subequations}
 The lowercase $s$, $c$ and $p$ denote concentrations of $S$, $C$ and 
$P$, respectively, and ``$\dot{\phantom{x}}$'' denotes differentiation with respect 
to time. The equation that models the time evolution of the enzyme concentration, $e$, has been eliminated via the enzyme conservation law, $e + c =e_0$, and note that adding together (\ref{eq:1a})--(\ref{eq:1c}) yields the substrate conservation law,
$s+c+p =s_0$. From the model equations (\ref{eq:1a})--(\ref{eq:1c}), we see that the MM 
reaction mechanism~(\ref{eq:mm1}) is parametrically controlled by the initial 
substrate concentration, $s_0$, the initial enzyme concentration, $e_0$, and the 
magnitudes of the rate constants $k_1$, $k_{-1}$ and $k_2$. 

Since the model equations (\ref{eq:1a})--(\ref{eq:1c}) are nonlinear, closed form 
solutions are intractable. However, it is well-established that whenever the 
quantity, $\varepsilon_{SS}$, is very small,
\begin{equation} \label{eq:segel-cond}
\varepsilon_{SS} = \cfrac{e_0}{K_M + s_0} \ll 1,
\end{equation}
equations (\ref{eq:1a})--(\ref{eq:1c}) can be \textit{approximated} by the system of
differential-algebraic equations:
\begin{equation}\label{eq:1red}
\dot{s} \sim -\cfrac{Vs}{K_M+s} \sim -\dot{p},\quad 
c \sim \cfrac{e_0s}{K_M+s} .
\end{equation}
Here the Michaelis constant is defined as $K_M \equiv (k_{-1}+k_2)/k_1$, 
and the limiting rate is $V\equiv k_2e_0$. Formally, the equations~(\ref{eq:1red}) 
are collectively referred to as the sQSSA, and condition (\ref{eq:segel-cond}) is known as the \textit{reactant-stationary-assumption} (RSA) \cite{Hanson:2008:RSA}. When the 
RSA holds, and the sQSSA (\ref{eq:1red}) is valid, the intermediate complex reaches its maximum value very quickly. By comparison, the 
conversion of $S$ to $P$ is slow, and the MM reaction 
mechanism~(\ref{eq:mm1}) is characterized (temporally) by the two disparate 
timescales, $t_C$ and $t_{\mathcal{D}}$, that respectively quantify the approximate amount of time it takes $c$ to become maximal, and for $s$ to significantly deplete. Thus, when the RSA holds, the reaction consists of a fast timescale, $t_C$, and a slow timescale, $t_{\mathcal{D}}$. Both timescales depend on the rate constants, as well as the initial substrate and enzyme concentrations:
\begin{subequations}\label{eq:scales}
\begin{align}
t_C &= \cfrac{1}{k_1(s_0 + K_M)},\\
t_{\mathcal{D}} &= \cfrac{K_M+s_0}{V}.
\end{align}
\end{subequations}

Segel~\cite{Segel1988} first derived $t_C$ and $t_{\mathcal{D}}$ 
using heuristic methods. To estimate $t_C$, Segel noted, no doubt from the earlier work of Heineken, Tsuchiya, and Aris~\cite{HEINEKEN196795}, that if $e_0 \ll s_0$, then $c$ should reach its maximum value very quickly, and there should be correspondingly very little substrate depletion during the initial accumulation of $c$. Consequently, one can assume $s\approx s_0$ during the initial accumulation of $c$. To obtain an equation that models the initial increase in $c$, we set $s=s_0$ in (\ref{eq:1b}), which gives rise to a linear ordinary differential equation for $c$:
\begin{equation}\label{eq:fcdot}
\dot{c} \sim k_1(e_0-c)s_0 - (k_{-1}+k_2)c.
\end{equation}
Since (\ref{eq:fcdot}) is linear, its solution is easily attainable, and is given by:
\begin{equation}\label{eq:fsol}
c(t) \sim \varepsilon_{SS} s_0(1-\displaystyle e^{\displaystyle -t/t_C}).
\end{equation}
Thus, $t_C$ is the \textit{characteristic} timescale of the initial fast transient, and 
it is valid as long as $s$ is approximately constant during the initial accumulation of $c$. The reader should be aware of the fact that (\ref{eq:fcdot}) is only valid for timescales on the order of $t_C$. Hence, the solution (\ref{eq:fsol}) is only an initial, temporary approximation to the kinetics of (\ref{eq:mm1}) when the RSA is valid. As such, it is often referred to as an \textit{initial layer}. 

Since  $t_C$ was derived under the presupposition of a negligible depletion 
of $s$ during the initial accumulation if $c$, how then, can we determine a dimensionless parameter 
that must be small in order to ensure the depletion of $s$ is in fact negligible 
over $t_C$? Segel~\cite{Segel1988} reasoned that since the maximum rate of depletion 
of $s$ is identically $-k_1e_0s_0$, it suffices to demand that
\begin{equation}\label{eq:dRSA}
\{\max|\dot{s}| \cdot t_C \ll s_0\} \equiv \{\varepsilon_{SS} \ll 1\}
\end{equation}
hold in order to justify the use of (\ref{eq:fcdot}). Thus, the RSA~(\ref{eq:segel-cond}) 
ensures that the approximation (\ref{eq:fcdot}) is valid. 

After $c$ reaches its maximum value, the QSS phase of the reaction ensues, and $c$ starts to slowly deplete. During the depletion phase $\dot{c}$ is very small 
but non-zero, and thus we say that $c$ evolves in a QSS or is ``slaved'' by $s$ since, according to the system (\ref{eq:1red}), the approximate concentration of $c$ is parametrically determined by the concentration of $s$. 
To estimate the depletion timescale, Segel~\cite{Segel1988} assumed that the 
sQSSA~(\ref{eq:1red}) was valid once $t_C \lesssim t$, and that the rate of depletion 
of substrate was approximately
\begin{equation}\label{eq:QSSA2}
    \dot{s} \sim -\cfrac{Vs}{K_M+s} \leq 0.
\end{equation}
Since $s\approx s_0$ for $t\lesssim t_C$,
it stands to reason that the sQSSA (\ref{eq:QSSA2}) can be supplied 
with the boundary condition $s(t_C)=s_0$, and thus the maximum rate of depletion 
in the QSS phase of the reaction is approximately
\begin{equation}\label{eq:QSSrate}
    \max |\dot{s}| \approx \cfrac{Vs_0}{K_M+s_0}.
\end{equation}
Finally, to calculate the depletion timescale, Segel~\cite{Segel1988} divided the total 
change in substrate over the course of the reaction, $\Delta s = s_0$, by the 
maximum rate of depletion in the QSS phase (\ref{eq:QSSrate}):
\begin{equation}\label{eq:TD}
t_{\mathcal{D}} = \cfrac{|\Delta s|}{\max|\dot{s}|} 
= \cfrac{s_0}{\cfrac{V s_0}{K_M+s_0}} = \cfrac{K_M+s_0}{V}.
\end{equation}

As noted in the introduction, the effective use of the singular perturbation method relies on separation of fast and slow timescales, and the question that naturally arises from the heuristic line of reasoning arises is: Does the RSA (\ref{eq:segel-cond}) ensure separation of timescales? That is, does it hold that $t_C \ll t_{\mathcal{D}}$ whenever $\varepsilon_{SS}\ll1$ ? The answer is yes: the RSA not only ensures that there is negligible formation 
of product (or depletion of substrate) for $t\lesssim t_C$, it also ensures 
that the timescales, $t_C$ and $t_{\mathcal{D}}$, are widely separated:
\begin{equation}
    \cfrac{t_C}{t_{\mathcal{D}}} \equiv \epsilon = \cfrac{\varepsilon_{SS}}{\bigg(1+\cfrac{k_{-1}}{k_2}\bigg)\bigg(1+\cfrac{s_0}{K_M}\bigg)}\leq \varepsilon_{SS}.
\end{equation}
Thus, the RSA (i.e., the condition that $\varepsilon_{SS} \ll 1$) \textit{induces} the 
separation of fast and slow timescales. Consequently, the RSA is understood to be 
a sufficient condition for the validity of the sQSSA (\ref{eq:1red}), and 
the separation of fast and slow timescales (i.e., $\epsilon \ll 1$) is a necessary 
(but not sufficient) condition for the sQSSA.

\begin{figure}[hbt!]
\centering
\includegraphics[width=10cm]{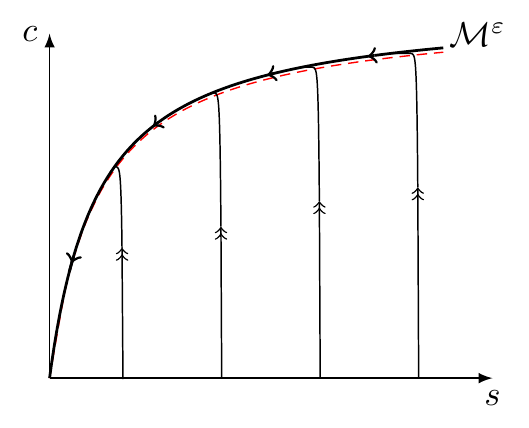}
\caption{\textbf{The RSA ensures the validity of the sQSSA: Fast/slow dynamics 
occur when $\varepsilon_{SS} \ll 1$ for the MM reaction mechanism~(\ref{eq:mm1})}. 
This figure is an illustration of the phase-plane dynamics when $\varepsilon_{SS}\ll1 $ and the sQSSA is valid; it serves to convey conceptual aspects of the phase-plane when the RSA is valid. The red dashed curve is the $c$-nullcline. The black curve the lies just above the $c$-nullcline is the invariant slow manifold. The thin lines with double arrows are illustrative of typical phase-plane trajectories when $\varepsilon_{SS}\ll 1$: trajectories starting on the $s$-axis approach the $c$-nullcline \textit{almost} vertically. Once the trajectory reaches the $c$-nullcline it closely follows $\mathcal{M}^{\varepsilon}$ during the QSS phase of the reaction. Double arrows indicate fast dynamics, and single arrows represent slow dynamics. The reader should bear in mind that, slightly contrary to what is illustrated, the tangent vector to the trajectory is perfectly horizontal when it intercepts the $c$-nullcline, since its vertical velocity component vanishes at the point of interception.}
\label{fig:sQSSA}
\end{figure}

When the RSA is valid, the trajectories in the $(s,c)$-plane have a very 
recognizable form: during fast transient phase, the trajectories move almost vertically 
towards the $c$-nullcline, given by the curve $\{(s,c)\in \mathbb{R}^2: c-e_0s/(K_M+s)=0\}$, on which $\dot{c}=0$.
The intermediate complex $c$ reaches its maximum value once the trajectory reaches 
the $c$-nullcline, at which point the QSS phase begins and the sQSSA is valid. After $c$ reaches 
its maximum value, the trajectory closely follows the $c$-nullcline towards the 
equilibrium point $(s,c)=(0,0)$ (see {\sc figure \ref{fig:sQSSA}}). 

The fact that the phase-plane trajectory \textit{follows} the $c$-nullcline on 
the slow timescale, and approaches the $c$-nullcline in vertical fashion over 
the fast timescale is the hallmark of the sQSSA. While the trajectory 
in {{\sc figure \ref{fig:sQSSA}}} appears as though it is \textit{on} the $c$-nullcline, it 
is important to note that it is not. In fact, the trajectory approaches an invariant manifold, $\mathcal{M}^{\varepsilon}$, that lies just above the $c$-nullcline. However, since $\mathcal{M}^{\varepsilon}$ and the 
$c$-nullcline are very close together when (\ref{eq:segel-cond}) holds, the 
$c$-nullcline is used to \textit{approximate} $\mathcal{M}^{\varepsilon}$ (again, see {\sc figure \ref{fig:sQSSA}} for a detailed 
explanation and illustration). The accuracy of using the $c$-nullcline as an approximation is assessed asymptotically through scaling and non-dimensionalization 
methods, and we review these methods in the subsection that follows.

\subsection{Asymptotic justification for the sQSSA: Scaling and
non-dimensionalization}
The identification of the small parameter~(\ref{eq:segel-cond}) as an 
appropriate condition for the validity of the standard QSSA is justified through 
scaling and non-dimensionalization. Introducing the dimensionless variables 
$\tau \equiv t/t_C, \bar{s} \equiv s/s_0$ and $\bar{c} \equiv c/\varepsilon_{SS} s_0$ 
yields
\begin{subequations}
\begin{align}
\cfrac{\text{d}\bar{s}}{\text{d}\tau} &= \varepsilon_{SS}[\mu\bar{c}\bar{s}-\bar{s} + \alpha\beta\bar{c}]\label{eq:sfast},\\
\cfrac{\text{d}\bar{c}}{\text{d}\tau} &=  \bar{s}-\mu\bar{c}\bar{s} - \beta\bar{c}\label{eq:cfast},
\end{align}
\end{subequations}
where $\beta \equiv 1/(1+\sigma)$, $\sigma \equiv s_0/K_M$, $\kappa \equiv k_{-1}/k_2$,
$\mu \equiv 1-\beta$ and $\alpha\equiv \kappa/(1+\kappa)$. Thus, we see from
(\ref{eq:sfast}) and (\ref{eq:cfast}) that if $\varepsilon_{SS}\ll 1$, then $s\approx s_0$ during the transient phase. 
Over the slow timescale, $T\equiv \varepsilon_{SS}\tau$, we obtain 
\begin{subequations}
\begin{align}
\cfrac{\text{d}\bar{s}}{\text{d}T} &= \mu \bar{c}\bar{s} -\bar{s} + \alpha\beta\bar{c}\label{eq:sslow},\\
\varepsilon_{SS}\cfrac{\text{d}\bar{c}}{\text{d}T} &=  \bar{s}-\mu\bar{c}\bar{s} - \beta\bar{c}\label{eq:cslow}.
\end{align}
\end{subequations}
Formally, 
the standard QSSA is an \textit{asymptotic} approximation of the dynamics on the 
slow timescale obtained by setting $\varepsilon_{SS} =0$ in (\ref{eq:cslow}) and solving for $c$ in terms of $s$. Thus, 
the standard QSSA is the zeroth-order approximation to (\ref{eq:sslow})--(\ref{eq:cslow}) 
and, when $\varepsilon_{SS} \ll 1$, we refer to $s$ as a \textit{slow} variable, 
and $c$ as a \textit{fast} variable.

In summary, we have two takeaways from this subsection:
\begin{enumerate}[label=(\roman*)]
    \item The RSA is synonymous with $\varepsilon_{SS} \ll 1$, and 
    guarantees that there is a negligible loss of substrate during the transient 
    (fast) phase of the reaction. When the RSA is valid, the QSSA system~(\ref{eq:1red}) 
    is an appropriate asymptotic approximation to (\ref{eq:1a})--(\ref{eq:1c}). 
    \item Timescale separation is synonymous with $\epsilon \equiv t_C/t_{\mathcal{D}}\ll 1$. 
\end{enumerate}
While (i) certainly implies (ii), the converse does not hold, which begs the 
question: What happens when $\epsilon \ll 1$ but $\varepsilon_{SS} \sim 1$? We discuss 
this special case in the subsection that follows.

\subsection{The \textit{extended} Quasi-Steady-State Approximation}
\begin{figure}[hbt!]
\centering
\includegraphics[width=8cm]{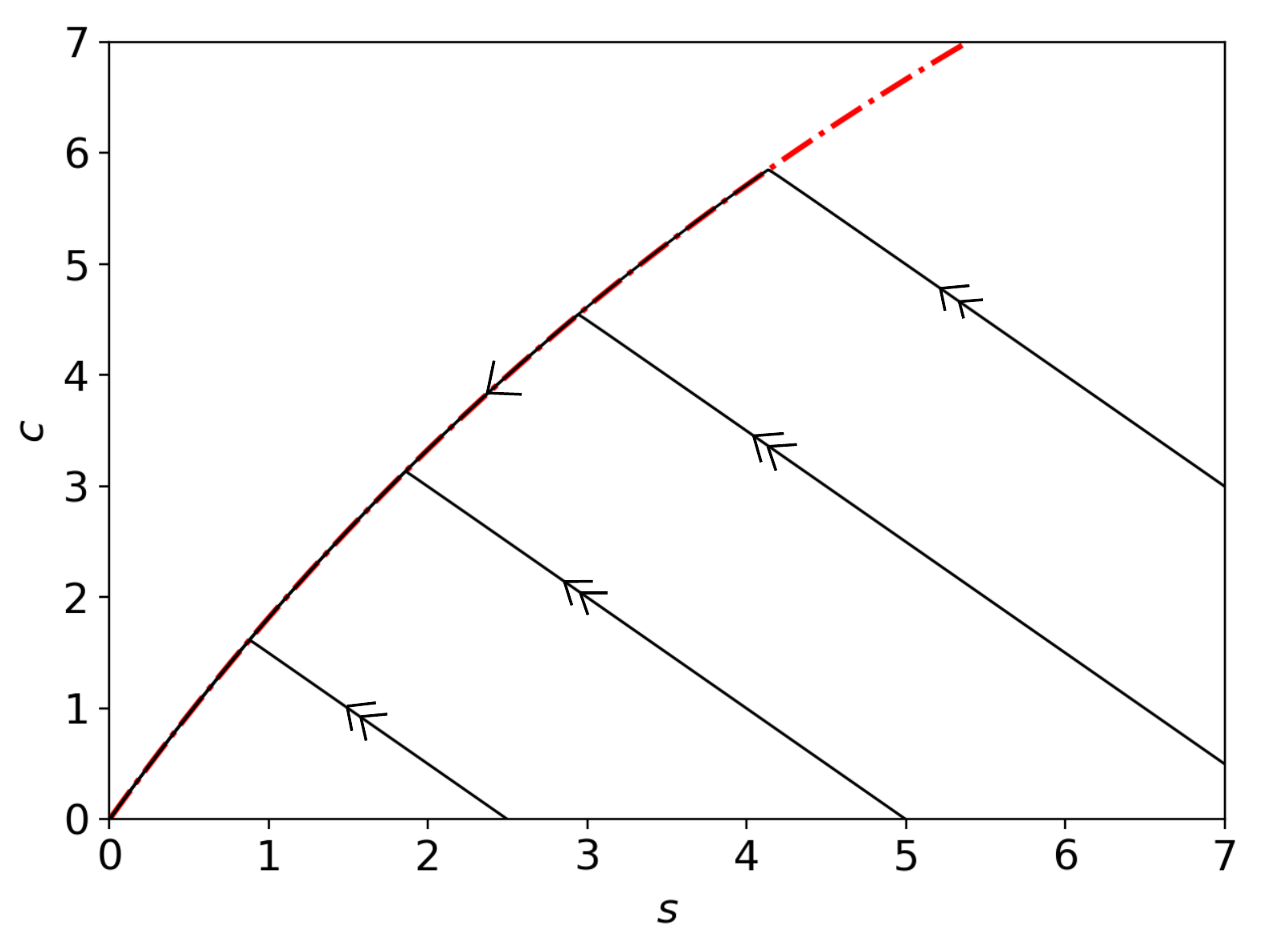}
\caption{\textbf{The extended QSSA: Fast/slow dynamics still occur when 
$\varepsilon_{SS} \sim 1$, but $\nu\ll 1$ for the MM reaction mechanism~(\ref{eq:mm1})}. 
The solid black curves are the numerical solutions to the mass action 
equations~(\ref{eq:1a})--(\ref{eq:1c}), and the dashed/dotted red curve 
is the ${c}$-nullcline. In this simulation, $k_1=k_{-1}=10$, $k_2=0.01$ 
with $\kappa \equiv k_{-1}/k_2 = 100$, $\nu = 1/101$, and  $\varepsilon_{SS} =1$
in arbitrary units for illustrative purposes. Note that the trajectory is not 
a vertical straight line during the initial phase of the reaction. As a 
consequence, the RSA does not hold~\cite{Segel1989,Hanson:2008:RSA}, and 
there is noticeable depletion of substrate during the transient phase. However, 
the trajectory still follows the $c$-nullcline after a brief transient and 
achieves the slow QSS dynamics. The double arrows indicate 
fast dynamics, and single arrows indicate slow dynamics.}
\label{fig:eQSSA}
\end{figure}

The RSA and timescale separation provide two different small parameters: $\varepsilon_{SS}$ 
and $\epsilon \equiv t_C/t_{\mathcal{D}}$. While $\varepsilon_{SS} \ll 1$ ensures 
that $\epsilon \ll 1$, the converse is not true. The obvious question is: Do QSS 
dynamics still prevail when $\varepsilon_{SS} \sim 1$ but $t_C/t_{\mathcal{D}}\ll 1$? 
Segel and Slemrod~\cite{Segel1989} proposed that the sQSSA~(\ref{eq:1red}) 
should still be valid at some point in the time course of the reaction as long as 
$t_C \ll t_{\mathcal{D}}$, even if $\varepsilon_{SS} \sim 1$. Rescaling the mass action 
equations (\ref{eq:1a})--(\ref{eq:1b}) with respect to $\bar{T}=t/t_{\mathcal{D}}$ yields
\begin{subequations}
\begin{align}
\nu\beta\cfrac{\text{d}\bar{s}}{\text{d}\bar{T}} &= \mu \bar{c}\bar{s} -\bar{s} + \alpha\beta\bar{c}\label{eq:sslowAA},\\
\nu\beta\varepsilon_{SS}\cfrac{\text{d}\bar{c}}{\text{d}\bar{T}} &=  \bar{s}-\mu\bar{c}\bar{s} - \beta\bar{c}\label{eq:cslowAA},
\end{align}
\end{subequations}
where again, $\nu = k_2/(k_{-1}+k_2)$, and $\beta = K_M/(K_M+s_0)$. If $\varepsilon_{SS} =1$, 
then setting $\nu = 0$ (which implies $\alpha=1$) in (\ref{eq:sslowAA})--(\ref{eq:cslowAA}) 
yields a (dimensional) manifold:
\begin{equation}\label{eq:snull}
    c = \cfrac{e_0s}{K_S+s},
\end{equation}
where $K_S = k_{-1}/k_1$ is the enzyme-substrate dissociation constant. Interestingly,
the manifold~(\ref{eq:snull}) is identical to the $s$-nullcline. However, 
since the $c$-nullcline and the $s$-nullcline are indistinguishable in the limiting 
case that corresponds to $k_2=0$, it stand to reason that trajectories should closely follow both the $s$- and 
$c$-nullcline when $k_2$ is incredibly small, but non-zero. Numerical simulations 
confirm that phase-plane trajectories closely follow the $c$-nullcline after a 
brief transient (see, {\sc Figure \ref{fig:eQSSA}}). 

The $c$-nullcline ``following'' that occurs when $\nu \ll 1$ and $\varepsilon_{SS}, \beta \sim 1$ is defined by Segel and Slemrod~\cite{Segel1989} as the \textit{extended} 
QSSA. One obvious
difference between the sQSSA and the extended QSSA is that there is a 
noticeable amount of substrate depletion when $\nu\ll 1$ but $\varepsilon_{SS} \sim 1$. 
Segel and Slemrod~\cite{Segel1989} proposed that the QSSA system~(\ref{eq:1red}) 
should still hold after the transient, but that the primary difference between 
the sQSSA and the extended QSSA was that one could not take $s=s_0$ at 
the onset of the QSS phase in the case of the extended QSSA. Thus, Segel and 
Slemrod~\cite{Segel1989} proposed that, when the extended QSSA is valid but the 
RSA is invalid, the approximation
\begin{equation}
    \dot{s} \sim -\cfrac{Vs}{K_M+s},
\end{equation}
must be supplied with a boundary condition ``$s_0^e$'' that is less than $s_0$. 
They employed a graphical method to estimate $s_0^e$ in their original 
manuscript~\cite{Segel1989}, and found that when $\beta=\varepsilon =1$, 
$s_0^e \approx (\sqrt{2}-1)s_0$. Thus, they concluded that the QSS phase of the reaction 
could be approximated with the initial value problem: 
\begin{equation}\label{eq:eQSSA}
    \dot{s} \sim -\cfrac{Vs}{K_M+s}, \quad s(0) = (\sqrt{2}-1)s_0.
\end{equation}
Moreover, they found that $s_0^e\to 0$ as $\beta$ increases, and determined 
that the parameter regime where $\nu \beta \ll 1$ but $\varepsilon_{SS} \sim 1$ 
results in a rather ``uninteresting'' domain of applicability for extended 
QSSA\footnote{ Segel and Slemrod~\cite{Segel1989} did not take into account 
product formation, and assumed that the complete depletion of substrate was 
synonymous with the completion of the reaction. Unfortunately, this assumption
is false. For more details, we invite the reader to consult~\cite{JTB}.}. Hence, 
the extended QSSA is formally defined to be the case when $\nu \ll 1$, but 
$\varepsilon_{SS}, \beta \sim 1$. Up to this point, we have two dynamical 
regimes of interest:
\begin{enumerate}[label=(\roman*)]
    \item The RSA, in which $\varepsilon_{SS} \ll 1$ and (\ref{eq:1red}) holds.
    \item The extended QSSA, in which $\nu \ll 1$, but 
    $\varepsilon_{SS}, \beta \sim 1$ and \textit{possibly} (\ref{eq:eQSSA}) holds.
\end{enumerate}

As a final remark of this subsection, we note that while the extended QSSA 
seems reasonable, there are a few technical deficiencies that engulf the validity of the approximation. First, setting $\nu =0$ 
results in the coalescence of the $s$- and $c$-nullclines, which implies that 
$s$ and $c$ are in some sense both in a QSS after the fast transient. Second, 
there is no scaling justification for (\ref{eq:eQSSA}). Consequently, other 
than the observation that the phase--plane trajectory follows the $c$-nullcline 
after the transient phase, there is really no rigorous justification 
for (\ref{eq:eQSSA}). We will return to these observations and discuss 
them in more detail in a later section.

\subsection{The \textit{reverse} Quasi-Steady-State Approximation}
\begin{figure}[hbt!]
\centering
\includegraphics[width=10cm]{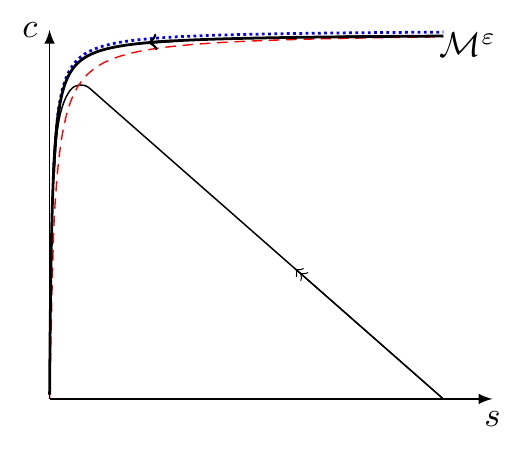}
\caption{\textbf{In the rQSSA, the substrate concentration is 
in a QSS phase 
in the MM reaction mechanism~(\ref{eq:mm1})}. This figure is an illustration of the rQSSA. The thick black curve is an illustration of the slow manifold; the dashed red curve in the $c$-nullcline, and the dotted blue curve is the $s$-nullcline. The thin black curve illustrates a typical trajectory when the rQSSA is valid. During the transient phase, the complex concentration rapidly reaches it maximum value, and the initial substrate concentration nearly vanishes. After the transient phase, the trajectory follows the slow manifold, $\mathcal{M}^{\varepsilon}$, during which time $s$ remains in a QSS phase. Double arrows represent fast dynamics, and single arrows represent slow dynamics.}
\label{fig:11}
\end{figure}

In contrast to the intermediate complex evolving in a QSS after the fast 
transient, the rQSSA occurs when the substrate evolves in QSS after the 
initial fast transient. Nguyen and Fraser~\cite{HNsF} observed that the 
rQSSA corresponds to trajectories closely following the $s$-nullcline 
instead of the $c$-nullcline in the $s$--$c$ phase-plane (see {\sc Figure~\ref{fig:11}}).  

To formulate an approximation to (\ref{eq:1a})--(\ref{eq:1c}), it is assumed that virtually all of the initial substrate 
will deplete during the initial transient phase, in which case
\begin{equation}\label{eq:rQ}
\cfrac{\text{d}c}{\text{d}t}+\cfrac{\text{d}p}{\text{d}t} =0
\end{equation}
\textit{approximately} holds during the QSS phase. It also assumed that there is negligible formation 
of product during the transient phase. Thus, in the rQSSA, the underlying 
assumption is that there is a preliminary transient phase in which nearly all 
of the substrate is depleted, but a negligible amount of product is generated. 
As a consequence of this assumption, (\ref{eq:rQ}) admits an approximate conservation law, $s_0\sim c+p$, that implies
\begin{equation}\label{eq:REA}
\cfrac{\text{d}p}{\text{d}t} = k_2(s_0-p)
\end{equation}
is the leading order approximation to $\dot{p}$ in the QSS phase of the 
reaction. The above expression governs the production formation under the 
rQSSA or \textit{rapid equilibrium 
approximation}~\cite{Roussel:1990:GSSA,Roussel:1991:ASS}. Equation (\ref{eq:REA}) 
is linear and admits a closed-form solution
\begin{equation}
p(t) = s_0(1-\displaystyle e^{\displaystyle -k_2t}),
\end{equation}
as well as a corresponding slow timescale: $\tilde{T}=k_2t$.

Is there a small parameter homologous to $\varepsilon_{SS}$ for the rQSSA that 
justifies the validity of (\ref{eq:REA})? Typically, less-than-rigorous heuristic arguments 
are employed to determine conditions for the validity of~(\ref{eq:REA}). By utilizing timescale and geometrical arguments, Segel and 
Slemrod~\cite{Segel1989} \textit{indirectly} proposed\footnote{ They 
directly proposed that $K_S/e_0\ll1$, but noted that $k_2\approx k_{-1}$, 
which implies $K_M/e_0 \ll 1$.}
\begin{equation}\label{eq:SSeps}
\varepsilon^* \equiv \cfrac{K_M}{e_0} = \eta^{-1} \ll 1
\end{equation}
as a necessary condition for the validity of the rQSSA. However, there are 
several problems with labeling (\ref{eq:SSeps}) as a \textit{sufficient} 
condition. First, as Schnell and Maini~\cite{SCHNELL2000} point out, 
$\varepsilon^* \ll 1$ does not ensure that $\text{d}s/\text{d}t$ vanishes 
after the initial transient, since $K_M \ll e_0 \ll s_0$ would seem to imply 
the validity of the QSSA, and Segel and Slemrod~\cite{Segel1989} only considered 
cases where $e_0 \gg s_0$ in their analysis of the rQSSA. Schnell and 
Maini~\cite{SCHNELL2000} suggested that there should be both a negligible loss 
of enzyme concentration and a maximal depletion of free substrate over the fast transient; 
they proposed   
\begin{equation}\label{eq:SM}
    \varepsilon_{SM} \equiv \eta^{-1} + \cfrac{s_0}{e_0} = \cfrac{1}{\varepsilon_{SS}} \ll 1
\end{equation}
as an appropriate qualifier for the validity of the rQSSA. The 
condition~(\ref{eq:SM}) suggests that since $e$ is in excess with respect 
to $s$, the majority of the substrate should vanish over the fast transient. Second, from a more technical perspective, 
Goeke et al.~\cite{Goeke2012} remark that $K_S/e_0 \ll 1$ is not sufficient 
for the validity of the rQSSA; we will say more about Goeke's observation in the section that follows. For now we note that, heuristically, if $1/k_2$ is a slow timescale, then it stands to reason that $k_2$ should also be small in some sense.

In short, 
most of the conditions in the literature \textit{basically} imply that 
$\varepsilon_{SS} \gg 1$ in order for the rQSSA to be valid. In this work, the question we want to address through our analysis is: How 
small should $k_{-1}$ and $k_2$ be in comparison to other variables to ensure 
the validity of the rQSSA (i.e., what is the equivalent of $\varepsilon_{SS}$ in 
this case)? Clearly, it must hold that $\varepsilon^* \ll 1$, but how does 
the initial concentration of substrate, $s_0$, impact the accuracy of the 
approximation (\ref{eq:REA})? This is still an open problem and, using 
non-scaling methods, we will close it in a later section. 

\subsection{The \textit{total} Quasi-Steady-State Approximation}
The tQSSA provides an approximation that is valid over a broader parameter domain. Utilizing the conservation law, $s_0=s+c+p$, and replacing 
$s=s_0-(c+p)$ in (\ref{eq:1b}) yields
\begin{subequations}
\begin{align}
\dot{c} &= k_1(e_0-c)(s_0-c-p) -k_1K_M c,\label{eq:tQSSA10}\\
\dot{p} &= k_2c\label{eq:tQSSA20}.
\end{align}
\end{subequations}
Setting the left hand side of (\ref{eq:tQSSA10}) to zero and solving for 
$c$ gives
\begin{equation}\label{eq:roots}
    c= \cfrac{1}{2}(e_0+K_M+s_0-p) \pm \cfrac{1}{2}\sqrt{(e_0+K_M+s_0-p)^2-4e_0(s_0-p)}.
\end{equation} 
Of the two possible roots given in (\ref{eq:roots}), only the ``$-$'' root 
is \textit{physical}. Taking
\begin{equation}\label{eq:mroot}
    c \approx \cfrac{1}{2}(e_0+K_M+s_0-p) - \cfrac{1}{2}\sqrt{(e_0+K_M+s_0-p)^2-4e_0(s_0-p)} \equiv c_0(p),
\end{equation}
and inserting it into (\ref{eq:mroot}) yields the tQSSA:
\begin{equation}\label{eq:tQ}
    \dot{p} \approx k_2 c_0(p).
\end{equation}
Thus, the complex is taken to be the fast variable that evolves in a QSS after a brief transient, and $p$, rather than $s$,  is the
corresponding slow variable. Due to the fact that $p$ is the designated 
slow variable, the tQSSA is valid over much larger domain than either 
the sQSSA or the rQSSA. As a result of the more recent work of
Borghans et al.~\cite{BORGHANS1996}, the term $s_0-p$ is usually replaced by $s_T=s+c$  in the literature (see, for example, \cite{SCHNELL2002, Pedersena2006, RAMITZAFRIRI2007737, Dingee, DellAcqua2012, Dhatt2013, KIM}), where $s_T$ is a lumped variable that denotes the ``total''
substrate. Consequently, the equations
\begin{subequations}
\begin{align}
\dot{c} &=  k_1(e_0-c)(s_T-c) -k_1K_M c,\\
\dot{s}_T &= k_2c,
\end{align}
\end{subequations}
are more commonly employed than (\ref{eq:tQSSA10})--(\ref{eq:tQSSA20}). 
Either way, as Pedersen \cite{PEDERSEN20081010} points out, the coordinate
systems $(s_T,c)$ and $(p,c)$ are merely ``two sides of the same coin''. We 
will employ the $(p,c)$ coordinate system, and avoid the lumped variable
``$s_T$'' in this paper as it is not a uniquely distinguishable physical
variable in experimental assays.

The validity of the tQSSA was first analyzed in the $(p,c)$ coordinate 
system by Laidler~\cite{Laidler:1955:TTK}, and later by Borghans et
al.~\cite{BORGHANS1996}. However, the accepted condition for the validity 
of the tQSSA was derived by Tzafriri~\cite{Tzafriri2003}. Briefly, Tzafriri
noted that if there is zero formation of product during the initial 
transient, then the formation of $c$ during the transient phase can be 
approximated by the Riccati equation
\begin{equation}\label{eq:RIC}
    \dot{c} =k_1(e_0-c)(s_0-c)-k_1K_Mc.
\end{equation}
The solution to (\ref{eq:RIC}) admits a natural \textit{fast} timescale:
\begin{equation}\label{eq:TCfast}
t_C^*=\cfrac{1}{k_1\sqrt{(e_0+K_M+s_0)^2-4e_0s_0}}.
\end{equation}
The corresponding \textit{slow} timescale, $t_P$, is calculated from 
the method of Segel~\cite{Segel1988}
\begin{equation}
    t_P = \cfrac{|\Delta p|}{\max |k_2 c_0(p)|} = \cfrac{2s_0}{(e_0+K_M+s_0) - \sqrt{(e_0+K_M+s_0)^2-4e_0s_0}},
\end{equation}
and the separation of $t_C^*$ and $t_P$ is the accepted condition 
for validity of the tQSSA:
\begin{equation}\label{eq:TZA}
    \varepsilon_T \equiv \cfrac{t_C^*}{t_P} = \cfrac{(e_0+K_M+s_0) - \sqrt{(e_0+K_M+s_0)^2-4e_0s_0}}{2s_0k_1\sqrt{(e_0+K_M+s_0)^2-4e_0s_0}} \ll 1.
\end{equation}
Several scaling analyses have been carried out that justify the asymptotic
validity of the tQSSA, and we refer the reader to Schnell and
Maini~\cite{SCHNELL2002} and Dell`Acqua and Bersani~\cite{DellAcqua2012} for detailed perturbation analyses, and Bersani
et al.~\cite{Bersani2015} for a thorough review of the topic. Since the 
tQSSA is valid whenever there is negligible formation of product during the fast
transient, the validity of the sQSSA, as well as the rQSSA and 
extended QSSA, implies the validity of the tQSSA. In this sense, the tQSSA
\textit{contains} both the standard QSSA and the rQSSA. Thus, 
$\varepsilon_{SS}, \nu \ll 1\implies \varepsilon_T\ll 1$, but $1\ll \varepsilon_{SS}$ also implies $\varepsilon_T \ll 1$ (again, see~\cite{BORGHANS1996}). 

It is worth noting that the fast and slow timescales derived by Tzafriri are valid whenever there is negligible product formation during the fast transient. Thus, they are valid in both the sQSSA and rQSSA regimes. Additionally, the upper bound on $c$ is easily obtained:
\begin{equation}
    \sup c \equiv \lambda = c_0(0) = \cfrac{1}{2}(e_0+K_M+s_0)-\cfrac{1}{2}\sqrt{(e_0+K_M+s_0)^2-4e_0s_0}.
\end{equation}
Consequently, we obtain universal timescales as well an upper bound on $c$ from the tQSSA formulation. At first glance the timescales $t_C^*$ and $t_P$, as well as the upper bound $\lambda$, appear quite complicated, and their utility in the context of scaling analysis appears limited. However, when we introduce non-scaling methods, we will show that in fact these ``ingredients" are quite useful. 

The tQSSA wraps up our introduction to each ``QSSA" (i.e., the sQSSA, rQSSA, and tQSSA) that is employed as a reduced model in enzyme kinetics. In summary, we have:
\begin{enumerate}[label=(\roman*)]
    \item  The RSA is valid when $\varepsilon_{SS} \ll 1$, and the corresponding QSSA is given by~(\ref{eq:1red}).
    \item  The extended QSSA is valid when when $\nu \ll 1$, but
    ($\varepsilon_{SS}\sim 1$ and $\beta \sim 1$). The corresponding QSSA is \textit{possibly} given by
    equation~(\ref{eq:eQSSA}).
    \item  The rQSSA is at least valid when $1\ll \varepsilon_{SS}$,
    and the corresponding QSSA is given by~(\ref{eq:REA}).
    \item  The tQSSA is valid whenever $\varepsilon_T \ll 1$, and the corresponding QSSA is given by~(\ref{eq:tQ}).
\end{enumerate}

\section{Applying the Quasi-Steady State Approximations with geometric singular 
perturbation theory}\label{sec:GSPT}
The work of Fenichel~\cite{GSPT0} consists of a group of theorems that warrant 
the existence of a slow invariant manifold, $\mathcal{M}$, in fast/slow 
systems of the form of (\ref{eq:1})--(\ref{eq:2}). In Section~\ref{sec:invariance}, we give a \textit{basic} introduction
to the results obtained by Fenichel and introduce relevant terminology for
our paper. For a more detailed introduction, we refer the reader 
to \cite{GSPT1,GSPT2} and \cite{kuehn2015multiple}. In Section~\ref{sec:GSPT},
we apply Geometric Singular Perturbation Theory (GSPT) to
(\ref{eq:1a})--(\ref{eq:1c}) under the assumption that 
the RSA is valid. In Section~\ref{sec:e-GSPT}, we apply GSPT to
(\ref{eq:1a})--(\ref{eq:1c}) in the parameter-space region where the 
extended QSSA is valid, and we show how GSPT can be used to prove that 
Segel and Slemrod's conclusion that the approximation given (\ref{eq:eQSSA}) 
is valid after the initial transient is incorrect.

\subsection{The invariance equation} \label{sec:invariance}
The sQSSA, rQSSA, and the tQSSA are, formally, zeroth-order approximations 
to the reaction dynamics on the corresponding slow timescale. However, 
there is much more to the story. Whenever the tQSSA is valid, trajectories 
are attracted to a \textit{slow} invariant manifold (SIM). The rigorous
study of SIMs is referred to as geometric singular perturbation theory 
(GSPT). Briefly, given a system of the form
\begin{subequations}
\begin{align}
    \dot{x}&=f(x,y;\varepsilon),\label{eq:A1} \\
    \varepsilon \dot{y}&=g(x,y;\varepsilon)\label{eq:B1},
\end{align}
\end{subequations}
where $x$ is a slow variable and $y$ is a fast variable, one can approximate 
the slow manifold by first assuming that the fast variable is expressible 
in terms of the slow variable as $y=h(x)$\footnote{The underlying assumption here is that the slow manifold can be expressed in terms of a global coordinate chart. More complicated manifolds, such as an n-dimensional torus, may require an atlas and a collection of transition maps.}, where $h$ is a function that 
is to be determined. Typically, the curve $y=h(x)$ is denoted as $\mathcal{M}^{\varepsilon}$, and we will use this notation from here 
on. Since the slow manifold is invariant, it must satisfy the differential
equation (\ref{eq:A1})--(\ref{eq:B1}). Hence, $y=h(x)$ must satisfy
\begin{equation}
    \cfrac{\text{d}h}{\text{d}x} \cdot f(x,h(x);\varepsilon) = \cfrac{1}{\varepsilon} g(x,h(x);\varepsilon),
\end{equation}
which is known as the \textit{invariance equation}. The power of the invariance equation is that, given a scaled equation
and an appropriate small parameter, it allows us to systematically determine the 
higher-order asymptotic approximations to the SIM. The function $h(x)$ admits a unique asymptotic expansion 
in terms of $\varepsilon$:
\begin{equation}
    h(x) = h_0(x) + \varepsilon h_1(x) + \varepsilon^2 h_2(x) + \mathcal{O}(\varepsilon^3),
\end{equation}
where $h_0(x)$ satisfies $g(y,h_0(x);0)=0$. Thus, when $g(x,y)$ is not
explicitly dependent on $\varepsilon$, the nullcline associated with 
the fast variable commonly serves as a zeroth-order approximation to 
the slow manifold, 
$\mathcal{M}^{\varepsilon}_0 = \{(x,y)\in \mathbb{R}^2 :y=h_0(x)\}$, and formally
the asymptotic approximation to the slow manifold can be expressed as:
\begin{equation}\label{eq:exp}
\mathcal{M}^{\displaystyle \varepsilon} = \mathcal{M}_0^{\displaystyle \varepsilon} + \varepsilon \mathcal{M}_1^{\displaystyle \varepsilon} + \varepsilon^2 \mathcal{M}_2^{\displaystyle \varepsilon}+...+ \mathcal{O}(\varepsilon^N)+...
\end{equation}
The terms $\varepsilon^0, \varepsilon^1, \varepsilon^2,..$ in (\ref{eq:exp}) are called \textit{gauge functions}, and the terms $\mathcal{M}_0^{\displaystyle \varepsilon}, \mathcal{M}_1^{\displaystyle \varepsilon},...$ are the coefficients that accompany the gauge functions. Generally speaking, the slow manifold is not necessarily \textit{unique}, meaning there can be slow manifold\textit{s}. However, they are separated by a distance that is $\mathcal{O}(e^{-K/\displaystyle \varepsilon})$ where $K$ is $\mathcal{O}(1)$~\cite{kuehn2015multiple}. 

In addition to the slow manifold, the work of Fenichel also guarantees the existence of \textit{fast fibers}, $\mathcal{F}^{\displaystyle \varepsilon,p}$, which can be utilized to approximate the dynamics during the transient stage as the trajectory rapidly approaches the slow manifold.\footnote{ In contrast 
to the critical manifold, the individual zeroth-order fast fibers to not perturb to individual invariant fibers when $\varepsilon$ is non-zero. However, the collective family of fast fibers,
$\displaystyle \bigcup_p \mathcal{F}^{\displaystyle \varepsilon,p}$, \textit{is} 
invariant in a certain way. For more details, we invite the reader to
consult \cite{Hek2010}.} The zeroth-order fiber,
$\mathcal{F}_0^{\displaystyle \varepsilon,p}$, which is the 
zeroth-order approximation to the fibers $\mathcal{F}^{\displaystyle \varepsilon,p}$, 
is the solution to the fast subsystem:
\begin{subequations}
\begin{align}
    x' &=0\label{eq:ff1},\\
    y' &= f(x,y;0)\label{eq:ff2}.
\end{align}
\end{subequations}
The fast subsystem contains a continuous branch of equilibrium points given by the curve $y=h_0(x)$, which is referred to in this context as the \textit{critical manifold}. Thus, each zeroth-order fiber links an initial condition $(x_0,y_0)$ to a unique \textit{base point} (equilibrium point) $p\in \mathcal{M}_0^{\displaystyle \varepsilon}$, 
and the superscript $p$ in $\mathcal{F}^{\displaystyle \varepsilon,p}$ is a reference 
to the unique base point, $p$. Since we are only interested in the 
zeroth-order approximation to $\mathcal{F}^{\displaystyle \varepsilon,p}$, we will 
henceforth drop the superscript $p$ for convenience, and use ``$\mathcal{F}_0^{\displaystyle \varepsilon}$'' to denote an \textit{individual} 
zeroth-order fiber. We will say more about fast fibers and the 
perspective of the slow manifold from the context of the fast 
subsystem when we introduce the idea of Tikhonov-Fenichel parameters 
in Section~\ref{sec:TK-parameters}. However, for now we will just 
associate the zeroth-order fast fiber as being the solution 
to (\ref{eq:ff1})--(\ref{eq:ff2}).

\subsection{Standard Quasi-Steady-State dynamics: Geometric singular perturbation 
theory when \texorpdfstring{$\varepsilon_{SS} \ll 1$}{Lg}.} \label{sec:sGSPT}
The RSA ensures negligible depletion of substrate during the fast 
transient, and \textit{suggests}\footnote{ This is, of course, a violation 
of the conservation law for substrate in the MM reaction 
mechanism~(\ref{eq:mm1}), but we will proceed as if it were permissible, and discuss this observation in more detail in {{section \ref{sec:nonscaled}.}}}
that the zeroth-order fast fiber, 
$\mathcal{F}_{0}^{\displaystyle \varepsilon_{SS}}$, is a straight line 
that connects the initial condition $(\bar{s},\bar{c})=(1,0)$ to the 
base point $(\bar{s},\bar{c})=(1,1)$:
\begin{equation}
\mathcal{F}_0^{\displaystyle \varepsilon_{SS}} = \{(\bar{s},\bar{c})\in 
\mathbb{R}^2 : \bar{s} = 1, 0\leq \bar{c} \leq 1\}.
\end{equation}

If $0 < \varepsilon_{SS} \ll 1$, then the time it takes the phase-plane
trajectory to reach a $\mathcal{O}(\varepsilon_{SS})$-neighborhood of 
the slow manifold is $\tau_{\varepsilon_{SS}} \sim |\ln \varepsilon_{SS}|$
or, in dimensional time $t \sim t_C |\ln \varepsilon_{SS}|$. Moreover, 
there exists an invariant manifold, $\mathcal{M}^{\displaystyle \varepsilon_{SS}}$, that is $\mathcal{O}(\varepsilon_{SS})$ from the 
critical manifold, and admits a unique expansion in powers of 
$\varepsilon_{SS}$ in the form of (\ref{eq:exp}),
where $\mathcal{M}_0^{\varepsilon_{SS}}$ is identically the $\bar{c}$-nullcline. These terms are straightforward 
to compute by constructing an asymptotic solution to the invariance 
equation:
\begin{equation}\label{eq:invariance}
\cfrac{1}{\displaystyle \varepsilon_{SS}}g(\mathcal{M}^{\displaystyle \varepsilon_{SS}},\bar{s}) = D_{\bar{s}}\mathcal{M}^{\displaystyle \varepsilon_{SS}}f(\mathcal{M}^{\displaystyle \varepsilon_{SS}},\bar{s}),
\end{equation}
and by inserting (\ref{eq:exp}) into (\ref{eq:invariance}),  it is easy 
to recover
\begin{equation}\label{eq:criticalnull}
    \mathcal{M}_0^{\varepsilon_{SS}} = \{(\bar{s},\bar{c})\in \mathbb{R}^2_+: \bar{c}-\bar{s}(1+\sigma)/(1+\bar{s}\sigma)=0\}.
\end{equation}
Finally, inserting the expression for $\mathcal{M}_0^{\varepsilon_{SS}}$ 
into (\ref{eq:sslow}) yields, in dimensional variables, the MM equation
\begin{equation}
\cfrac{{\rm d}s}{{\rm d}t} \sim -\cfrac{Vs}{K_M+s},
\end{equation}
which approximates the long-time dynamics of the mass action equations when the RSA is valid.

\subsection{Extended Quasi-Steady-State dynamics: Geometric singular 
perturbation theory when \texorpdfstring{$\nu\ll 1$, $\beta\sim 1$}{Lg} 
and \texorpdfstring{$\varepsilon\sim 1$}{Lg}.} \label{sec:e-GSPT}
Now let us consider the case of the extended QSSA. The problem can be
reformulated in the $(p,c)$ coordinate system so that it is of the form (\ref{eq:A1})--(\ref{eq:B1}) (see, for example, \cite{HEINEKEN196795,ROUSSEL2019108274}).
However, we will work out the solution in the $(s,c)$ coordinate system. 
In this case, we can guess the asymptotic expansion of the slow 
manifold in this coordinate system, and assume it takes the form
\begin{equation}\label{eq:cnullEXP}
c = \mathcal{M}^{\displaystyle \tilde{\nu}} = \cfrac{\bar{s}}{\alpha\beta + \mu\bar{s}} + \nu \mathcal{M}^{\displaystyle \tilde{\nu}}_1 + \mathcal{O}(\displaystyle \tilde{\nu}^2), \quad \tilde{\nu} \equiv 1/\kappa = k_2/k_{-1},
\end{equation}
where the curve
\begin{equation}\label{eq:snull-zeroth}
    \bar{c}=\cfrac{\bar{s}}{\alpha\beta + \mu\bar{s}},
\end{equation}
is the dimensionless form of the $s$-nullcline (\ref{eq:snull}). Taking 
just the zeroth-order expansion of (\ref{eq:cnullEXP}) (i.e., (\ref{eq:snull-zeroth})) and inserting it into equation (\ref{eq:sslowAA}) 
while zeroing all terms of order $\nu$ or higher yields
\begin{equation}
\cfrac{\text{d}\bar{s}}{\text{d}T} = 0.
\end{equation}
Since the $\bar{c}$- and $\bar{s}$-nullclines coalesce when $k_2 =0$, the zeroth-order approximation to the dynamics on the slow manifold 
results in the flow being infinitely slow (i.e., an equilibrium solution). Consequently, we must go to first-order 
in $\tilde{\nu}$ to determine the QSS approximation to the dynamics 
on the slow manifold. Solving for 
$\mathcal{M}^{\displaystyle \tilde{\nu}}_1$ by approximating the 
solution to the invariance equation\footnote{$D_{\bar{s}}$ denotes
differentiation with respect to $\bar{s}$.}
\begin{equation}
g(\mathcal{M}^{\displaystyle \tilde{\nu}},\bar{s}) = D_{\bar{s}}\mathcal{M}^{\displaystyle \tilde{\nu}}f(\mathcal{M}^{\displaystyle \tilde{\nu}},\bar{s})
\end{equation}
in powers of $\tilde{\nu}$, we recover
\begin{equation}\label{eq:skappa}
\cfrac{\text{d}\bar{s}}{\text{d}T} = \tilde{\nu}\cfrac{((\bar{s}-1)\theta +1)(\theta -1)\bar{s}}{(\bar{s}-1)^2\theta^2+(2\bar{s}-2\varepsilon -1)\theta + \varepsilon +1} + \mathcal{O}(\tilde{\nu}^2), \quad \theta = \cfrac{\sigma}{\alpha+\sigma} < 1,
\end{equation}
by inserting $\mathcal{M}^{\displaystyle \tilde{\nu}}_1$ into (\ref{eq:sslowAA}), which is the leading-order non-trivial solution 
to the dynamics on the slow manifold. Thus, the QSS
approximation for $s$~(\ref{eq:skappa}) is different when $\nu \ll 1$ versus 
when $\varepsilon_{SS} \ll 1$. In dimensional variables, 
(\ref{eq:skappa}) translates to,
\begin{equation}\label{eq:kappa}
\cfrac{\text{d}s}{\text{d}t}\sim-\cfrac{Vs(s+K_S)}{e_0K_S + (s+K_S)^2}, \quad K_S \equiv k_{-1}/k_1.
\end{equation}
The result (\ref{eq:kappa}) is not new, and the approximation that can be found in several papers \cite{Boulier2009,GoekeBook,GOEKE20171PhysicaD,Goeke2012,ROUSSEL2019108274}. The point is that by utilizing the 
invariance equation, we have shown that in fact Segel and Slemrod's 
extended QSS approximation~(\ref{eq:eQSSA}) is \textit{invalid}. The difference between the QSS approximations for $s$ when $\nu \ll 1$ versus when $\varepsilon_{SS} \ll 1$ has a nice history, and Roussel \cite{ROUSSEL2019108274} recently published a thorough explanation of the different QSS approximations through the lens of Tikhonov's Theorem and geometric singular perturbation theory. 

\begin{figure}[hbt!]
\centering
\includegraphics[width=8cm]{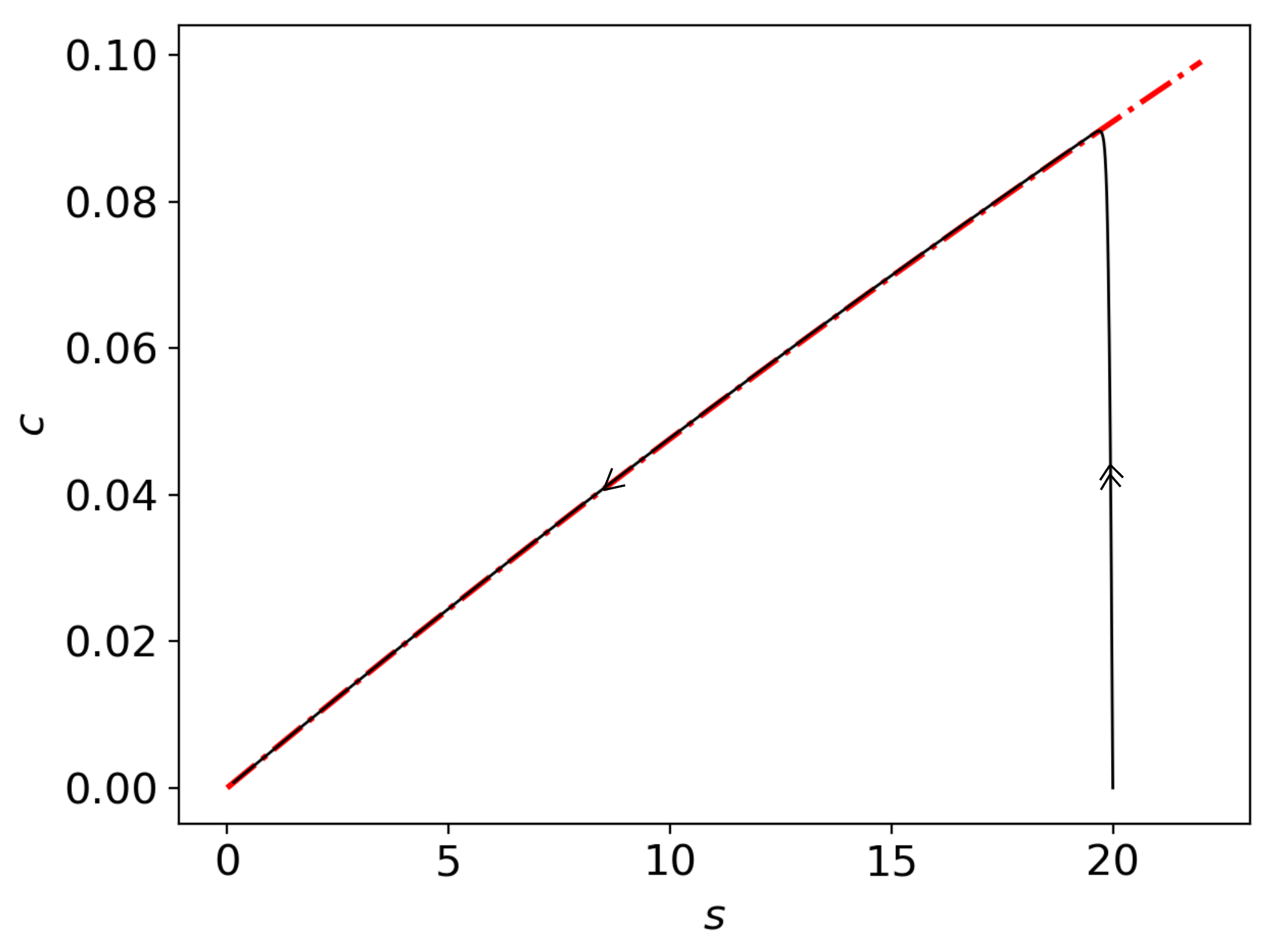}
\includegraphics[width=8cm]{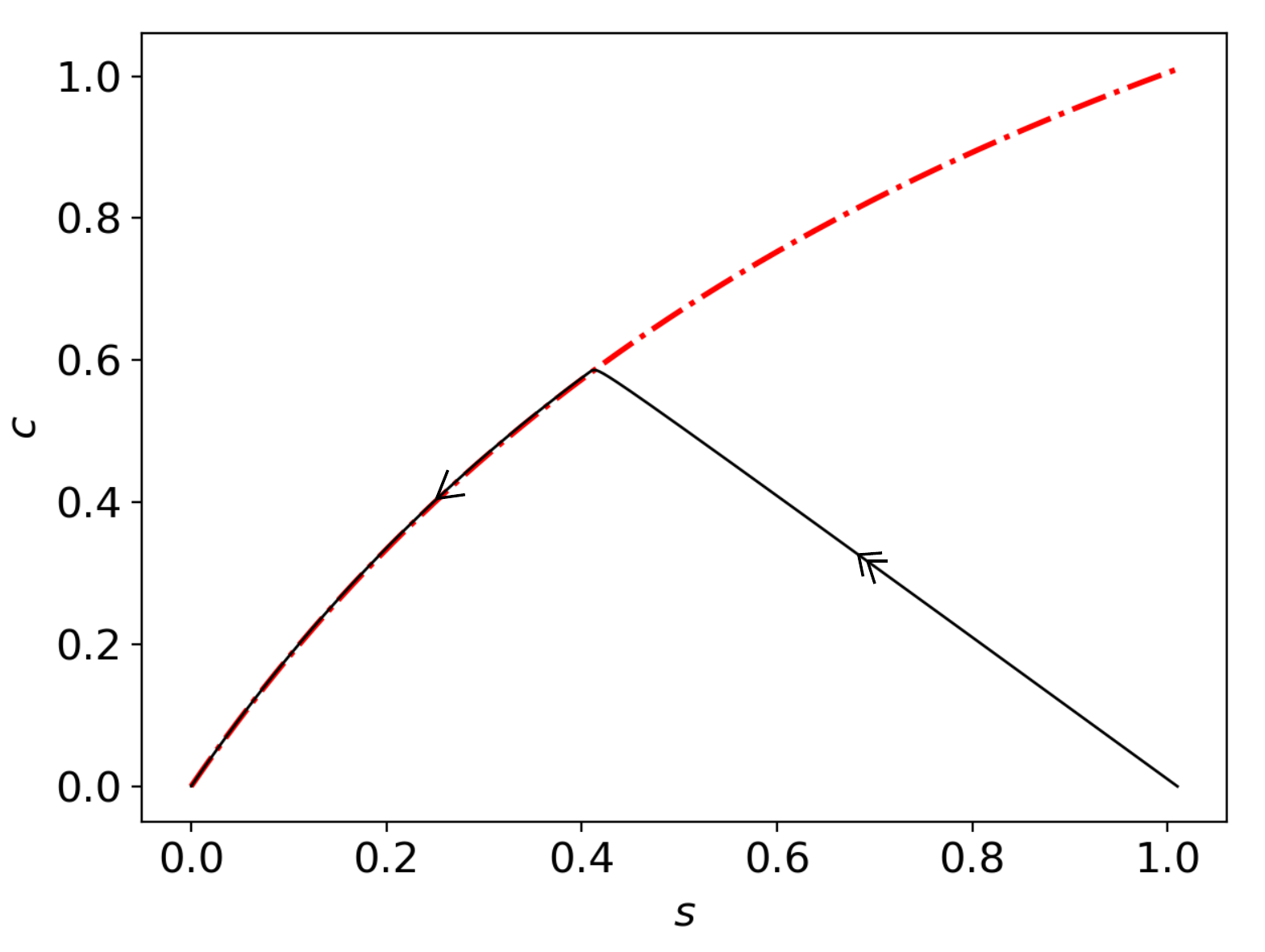}
\caption{\textbf{The extended QSSA equation~(\ref{eq:eQSSA}) proposed 
by Segel and Slemrod~\cite{Segel1989} is invalid when $\nu\ll 1$ but
$\varepsilon_{SS}=\sigma=1$}.
In both panels, the solid black curve is the numerical solutions to the mass action 
equations~(\ref{eq:1a})--(\ref{eq:1c}), and the dashed/dotted red curve 
is the $c$-nullcline. {\sc{left}}: In this panel, $s_0=20, e_0=1,k_2=k_{-1}=10, k_1=0.1$, and the sQSSA is valid, and the trajectory moves almost vertically towards the $c$-nullcline, which is the hallmark of the sQSSA in the $(s,c)$ phase-plane. {\sc{right}}: This panel is illustrates the path of a trajectory when $k_2 \ll k_{-1}$ but $\varepsilon_{SS} =1$. The parameters utilized in the numerical simulation are:
$k_1=k_{-1}=1,k_2=0.01$ with $\kappa = 100$ and $\sigma =1$. Notice the
trajectory is no longer \textit{vertical} during the initial phase of the reaction 
as there is substantial loss of substrate during the fast transient. Once 
the trajectory reaches the $c$-nullcline, the approximation 
given by (\ref{eq:kappa}) is valid, and can be equipped with the initial
condition $(\sqrt{2}-1)s_0$ if $\varepsilon_{SS} = \sigma =1$. Although the sQSSA \textit{looks} valid in both panels based on the proximity of the trajectory to the $c$-nullcline, GSPT indicates it is only valid in the left panel. In this figure, 
double arrows represent fast dynamics, and single arrows represent slow
dynamics. All units are arbitrary.}
\label{fig:21}
\end{figure}

In addition to the change in leading-order dynamics on the slow manifold 
when $\nu \ll 1$ instead of $\varepsilon_{SS}\ll 1$, the zeroth-order
approximation to the fast fiber, $\mathcal{F}_0^{\nu}$, will also change. One major difference between the transient dynamics of the sQSSA and the extended QSSA is that we can no longer take $s\approx s_0$ during the trajectory's initial approach to the slow manifold. Hence, we must determine an appropriate starting point in order to ``match" the inner and outer solutions. Locating an appropriate starting position on a slow manifold is not necessarily a trivial exercise. However, mathematicians have developed methods for such a situation (in~\cite{Roberts89}, A. J. Roberts develops a clever method for estimating starting positions on centre manifolds). Luckily, in our case, the conservation of total substrate, $s_0=s+c+p$, can be utilized to estimate a starting point.
Over the $\tau$-timescale, the total system is given by:
\begin{subequations}
\begin{align}
\cfrac{\text{d}\bar{s}}{\text{d}\tau} &= \varepsilon_{SS}[\mu\bar{c}\bar{s}-\bar{s} + \alpha\beta\bar{c}]\label{eq:sfasttt},\\
\cfrac{\text{d}\bar{c}}{\text{d}\tau} &=  \bar{s}-\mu\bar{c}\bar{s} - \beta\bar{c}\label{eq:cfasttt},\\
\cfrac{\text{d}\bar{p}}{\text{d}\tau} &= \epsilon \bar{c}\label{eq:pfasttt}.
\end{align}
\end{subequations}
Since $\nu \to 0$ implies $\epsilon \to 0$, we set $\epsilon =0$ and $\alpha = 1$ in (\ref{eq:pfasttt}) and (\ref{eq:sfasttt}), respectively. Setting the right hand side of (\ref{eq:pfasttt}) to zero implies that the \textit{total} substrate, 
``$\bar{s}_T = \bar{s} + \bar{c},$'' is conserved over the transient timescale when $\nu =0$,
\begin{equation}\label{eq:total}
\bar{s} + \bar{c} = 1,
\end{equation}
and insertion of (\ref{eq:total}) into the scaled equations yields a Riccati 
equation for $\bar{c}$:
\begin{equation}\label{eq:riccati}
\cfrac{\text{d}\bar{c}}{\text{d}\tau} = 1-2\bar{c} + \mu\bar{c}^2.
\end{equation}
Two observations can be made. First,
the correction to the duration of the fast transient when $\nu \ll 1$ (as opposed to $\varepsilon_{SS}$) can be found by 
solving the Riccati equation (\ref{eq:riccati}). In dimensional form, the solution to (\ref{eq:riccati}) yields (\ref{eq:TCfast}).
Again, this timescale was originally introduced by 
Tzafriri~\cite{Tzafriri2003}, and can be taken to be characteristic of 
the fast transient whenever there is minimal product formation during the time it takes $c$ to reach its maximum value (i.e., even when $\varepsilon_{SS} \ll 1$).\footnote{The usage of ``characteristic" is slightly abused in this context since the Riccati equation is nonlinear.} Second, the 
base point $(\bar{c}^*, \bar{s}^*)$ of the zeroth-order fast fiber is given by\footnote{ The zeroth-order fast fiber is not a vertical line 
in the $(s,c)$ phase--plane. However, it is a vertical line in the 
$(p,c)$ phase--plane, since $p$ is \textit{the} slow variable in this 
case and not $s$.}
\begin{equation}
\bar{c}^* = \cfrac{2-\sqrt{4(1-\mu)}}{2 \mu}, \quad \bar{s}^* = 1-\bar{c}^*.
\end{equation}
Consequently, if $\varepsilon_{SS} =\sigma =1 $ then, after the fast 
transient, we expect $\bar{c} \approx 2-\sqrt{2}$ and $\bar{s} \approx \sqrt{2}-1$ 
(see {\sc figure \ref{fig:21}}). In dimensional variables, this corresponds 
to $c \approx (2-\sqrt{2})\,s_0$ and $s = (2-\sqrt{2})\,s_0$, and we recover the 
original estimate given by Segel and Slemrod~\cite{Segel1989} (see 
{\sc figure \ref{fig:21}}).

GSPT, or really just the invariance equation, provides a powerful 
framework from which QSS approximations can be carefully derived. 
Moreover, it can help to rule out erroneous conclusions about
specific QSS approximations such as (\ref{eq:eQSSA}). Note that the 
utility of the invariance equation still relies on direct knowledge 
of a suitable small parameter in order to construct an asymptotic expansion 
to the slow manifold. However, the appropriate small parameter that warrants the validity of the rQSSA is still an open problem. This raises 
the question: are there suitable non-scaling 
methods that can be employed to determine the validity of the rQSSA? We introduce two such methods in the section that follows.

\section{Applying the Quasi-Steady-State Approximations to the 
Michaelis–Menten reaction mechanism: Non-scaled approaches}\label{sec:nonscaled}
In this section, we introduce non-scaling approaches to finding 
parameters or combinations of parameters that, when made very small, 
justify the validity of the sQSSA, the extended QSSA or the 
rQSSA. We also remark on the \textit{origins} of slow manifolds, as this will be critical to uncovering sufficient conditions for the validity of the rQSSA in Section \ref{sec:rQSSA}.

\subsection{Geometric singular perturbation theory: The fast subsystem}
\begin{figure}[hbt!]
\centering
\includegraphics[width=10cm]{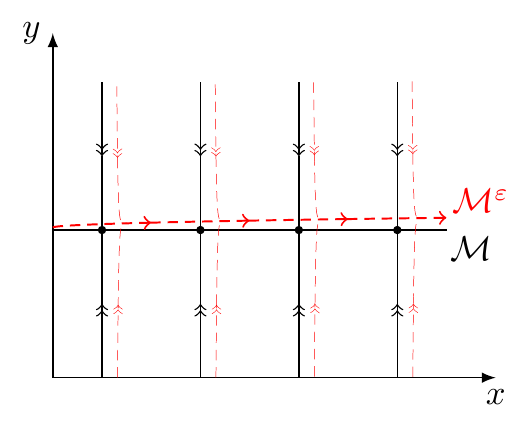}
\caption{\textbf{Illustration of a normally hyperbolic and invariant 
critical manifold, $\mathcal{M}$}. The black lines with arrows and fixed points represent trajectories within the unperturbed vector field field; the red dashed lines represent trajectories in the perturbed vector field. The straight black lines with 
double arrowheads represent trajectories moving towards the fixed 
points (filled black circles) that lie on the normally hyperbolic, 
attracting and invariant manifold $\mathcal{M}$, depicted by the 
thick black line. The flow on $\mathcal{M}$ is trivial. However, 
once $\varepsilon$ is non-zero and the vector field is smoothly 
perturbed, an invariant manifold $\mathcal{M}^{\varepsilon}$ 
emerges (red), on which the flow is slow, but no longer trivial. Trajectories that start off of $\mathcal{M}^{\varepsilon}$ quickly approach it, then follow it tangentially. Obviously, 
the fixed points drawn on this manifold are isolated, so the reader 
should bear in mind that $\mathcal{M}$ is filled with non-isolated 
fixed points.}
\label{fig:2121}
\end{figure}

Let us consider once more a fast/slow system of the form
\begin{subequations}
\begin{align}
x' &= \varepsilon f(x,y;\varepsilon)\label{eq:pert1},\\
y' &= g(x,y; \varepsilon)\label{eq:pert2},
\end{align}
\end{subequations}
which we will refer to as being in \textit{Tikhonov standard form.} 
The associated fast subsystem,
\begin{subequations}
\begin{align}
x' &= 0\label{eq:fastsub1A},\\
y' &= g(x,y;0)\label{eq:fastsub2A},
\end{align}
\end{subequations}
contains a branch of fixed points $y=h_0(x)$, where $g(x,h_0(x);0)=0$. 
Again, each fixed point in the set $y=h_0(x)$ is called a base point, and 
the complete set $y=h_0(x)$, denoted hereby as $\mathcal{M}$, is said to be
normally-hyperbolic provided
\begin{equation}\label{eq:hyper}
    \mathfrak{Re}\bigg( \cfrac{\partial g(x,y;0)}{\partial y}\bigg|_{y=h_0(x)}\bigg) \neq 0.
\end{equation}
In this context, normal hyperbolicity is partially\footnote{ Normally hyperbolicity 
provides information about the dominant directions of the flow near 
the critical manifold.} a statement about the linear stability along directions that are transverse to the critical manifold. As previously 
stated, the zeroth-order fast fibers are solutions to 
(\ref{eq:pert1})--(\ref{eq:pert2}), but they are also individual 
manifolds themselves. To see this, note that if (\ref{eq:hyper}) 
is less than zero, then $\mathcal{M}$ is attracting, and trajectories 
that start \textit{close} to $\mathcal{M}$ move closer to $\mathcal{M}$  as $t\to\infty$ (see {\sc Figure \ref{fig:2121}}). 
In contrast, if (\ref{eq:hyper}) is positive, then $\mathcal{M}$ 
is repelling, and trajectories that start close to $\mathcal{M}$ 
move away from $\mathcal{M}$ as $t\to \infty$. Thus, when the critical manifold is normally hyperbolic and attracting, 
the zeroth-order fast fiber is the stable manifold corresponding 
to the base point ``$p$'' of $\mathcal{M}$, and the stable manifold 
of $\mathcal{M}$ is comprised of (foliated by) the entire union of
zeroth-order fast fibers.  
Why are we interested in manifolds that are comprised of fixed 
points? Fenichel's Theorems describe what happens to normally 
hyperbolic critical manifolds of fixed points once the vector field is perturbed. 
Formally, we introduce the following theorem due to Fenichel~\cite{GSPT0},
which can be found in \cite{kuehn2015multiple}:

\textbf{Theorem 1:}
Let $F$ be a $C^r$ vector field on $\mathbb{R}^n$ with $r\geq 1$. 
Let $\mathcal{M}$ be a compact and connected $C^r$ manifold embedded 
in $\mathbb{R}^n$. Suppose that $\mathcal{M}$ is normally hyperbolic 
and invariant under the flow of $F$. Then, given any $C^r$ vector 
field $F^{\varepsilon}$ that is sufficiently $C^1$--close to $F$, 
there exists a $C^r$--manifold $\mathcal{M}^{\varepsilon}$ that is 
invariant under the flow of $F^{\varepsilon}$ and diffeomorphic 
to $\mathcal{M}$.

The notation ``$C^r$'' means the vector field has at least $r$ 
derivatives, all of which are continuous. The terms ``compact'' 
and ``connected'' are referring to specific topological properties 
that $\mathcal{M}$ must be equipped with. Connected means that 
$\mathcal{M}$ cannot be expressed as the union of two disjoint 
and non-empty sets, and compact means every open cover has a finite sub-cover. 
The term $C^1$-close is a statement about the distance between 
$F$ and $F^{\varepsilon}$:

\textbf{Definition:}
Let $F$ and $F^{\varepsilon}$ be two $C^1$ vectors fields on 
$\mathbb{R}^n$, and $\Lambda$ be a compact set. 
$F$ is said to be $\vartheta$-close to $F^{\varepsilon}$ on 
$\Lambda$ if:
\begin{subequations}
\begin{align}
    \sup_{\boldsymbol{x}\in \Lambda} &|| F-F^{\varepsilon}|| \leq \vartheta,\\
    \sup_{\boldsymbol{x}\in \Lambda} &|| DF - DF^{\varepsilon}|| \leq \vartheta,
\end{align}
\end{subequations}
where $D$ denotes differentiation. {\bf Theorem 1} is one of several 
theorems that provides the foundation for GSPT (see \cite{Hek2010} for a more technical survey). Later theorems address compact manifolds with a boundary, in which case invariance is replaced with \textit{local} invariance. In a nutshell, one 
starts with a vector field that contains a slow, compact critical manifold 
that consists of fixed points. Fenichel's theorems tells us that if a compact subset of the 
critical manifold is normally hyperbolic and attracts 
nearby trajectories, then perturbing the vector field (i.e., 
allowing $\varepsilon$ in (\ref{eq:fastsub1A})--(\ref{eq:fastsub2A}) 
to be greater than zero) in an appropriate way results in an 
attracting slow manifold. 

\begin{figure}[hbt!]
\centering
\includegraphics[width=8cm]{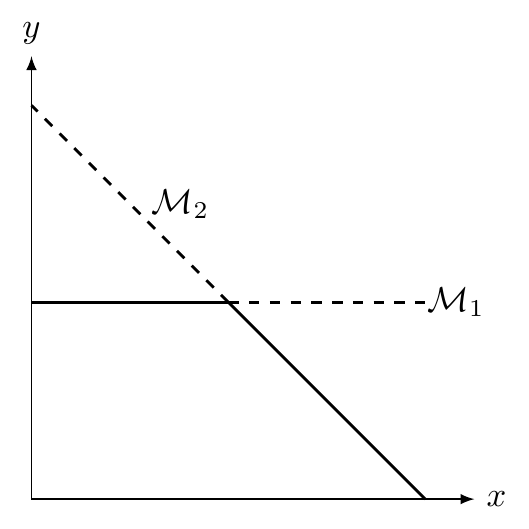}
\caption{\textbf{The transcritical bifurcation point corresponds 
to a loss of normal hyperbolicity and an exchange of stability}. Here we illustrate a transcritical
bifurcation. Thick dashed lines corresponds to unstable fixed points 
in which (\ref{eq:hyper}) is greater than zero. Thick solid lines
represent stable branches of fixed points for which (\ref{eq:hyper}) 
is negative. Two manifolds, $\mathcal{M}_1$ and $\mathcal{M}_2$, cross 
at a singular point. At the crossing point (i.e, the point of 
intersection), normal hyperbolicity is lost, and (\ref{eq:hyper})
fails to hold. Each manifold undergoes a change in stability at 
the point of intersection (i.e., the singular point delimits the 
point at which the stability of each of the manifolds changes from 
attractive to repulsive or vice-versa).}
\label{fig:ABC}
\end{figure}

Let us now discuss a specific point,  
``$p\in \mathcal{M}$'', where (\ref{eq:hyper}) fails and
\begin{equation}\label{eq:hyper1}
     \cfrac{\partial g(x,y;0)}{\partial y}\bigg|_{p\in\mathcal{M}} = 0.
\end{equation}
A point that satisfies (\ref{eq:hyper}) is called a \textit{singular point}. Singular points can indicate where a change in the stability of 
the critical manifold $\mathcal{M}$ may occur, and can correspond 
to dynamic bifurcation points. In particular, a transcritical 
bifurcation occurs at a singular point where two normally hyperbolic
submanifolds cross and undergo an exchange of stability (see 
{\sc Figure~\ref{fig:ABC}}). This specific bifurcation will 
be of interest when we analyze the rQSSA. 

\subsection{Tikhonov-Fenichel parameters} \label{sec:TK-parameters}
We now want to take the idea of a critical manifold that is composed 
of fixed points one step further. Revisiting the mass action equations
(\ref{eq:1a})--(\ref{eq:1c}), we can ask the following question: Are 
there specific dimensional parameters that, when identically zero, 
result in the formation of a normally hyperbolic invariant manifold 
that is composed of fixed points? If so, then it seems reasonable to 
assume that these parameters should in some sense be small whenever 
a specific QSS reduction is valid. For example, let us first consider 
the rQSSA in the $(s,p)$ phase-plane, for which the mass action 
equations are:
\begin{subequations}
\begin{align}
\dot{s} &= -k_1(e_0-(s_0-s-p))s + k_{-1}(s_0-s-p),\label{eq:sp1}\\
\dot{p} &= k_2(s_0-s-p)\label{eq:sp2}.
\end{align}
\end{subequations}
We know that the rQSSA implies that $s\approx 0$ during the QSS phase 
of the reaction. Define $k_2 = \varepsilon_2 k_2^{\star}$ and
$k_{-1}=\varepsilon_1 k_{-1}^{\star}$, where $\varepsilon_1$ and $\varepsilon_2$ are 
dimensionless and small, $0< \varepsilon_1,\varepsilon_2 \ll 1$, and 
$k_{-1}^{\star},k_2^{\star}$ are identically one with the same 
units as $k_{-1}$ and $k_2$  respectively. 
Expressing (\ref{eq:sp1})--(\ref{eq:sp2}) as
\begin{subequations}
\begin{align}
\dot{s} &= -k_1(e_0-(s_0-s-p))s + \varepsilon_1 k_{-1}^{\star}(s_0-s-p)\label{eq:sp12},\\
\dot{p} &= \varepsilon_2 k_2^{\star}(s_0-s-p)\label{eq:sp22},
\end{align}
\end{subequations}
it is clear that the curve $s=0$ is filled with equilibrium points 
when $\varepsilon_1=\varepsilon_2 =0$ and $e_0>s_0$. Furthermore, it is 
straightforward to show that the curve $s=0$ is normally 
hyperbolic and attracting when $s_0 < e_0$, since
\begin{equation}
    \cfrac{\partial}{\partial s} \bigg[ -k_1(e_0-(s_0-s-p))s\bigg]\bigg|_{s=0} < 0.
\end{equation}
Consequently, we anticipate
\begin{equation}
    \dot{p} \sim k_2(s_0-p)\label{eq:rQAAAA}
\end{equation}
during the slow QSS phase of the reaction when $0< \varepsilon_1,\varepsilon_2 \ll 1$. The takeaway message
from this example is two-fold. First, without scaling and
non-dimensionalization, we have predicted the validity of the rQSSA 
by choosing dimensional parameters (i.e., $k_{-1}$ and $k_2$) that, 
upon vanishing, give rise to an invariant manifold of equilibrium 
points. This is consistent with the notion of a normally hyperbolic 
critical manifold from fast/slow analysis and GSPT. Second, any 
dimensionless parameter homologous to $\varepsilon_{SS}$ that 
acts as a qualifier for the validity of the rQSSA should vanish 
whenever $k_{-1}$ \textit{and} $k_{2}$ vanish. In planar systems, a parameter that, upon vanishing, gives rise to a manifold of fixed points that is normally hyperbolic, is called a
\textit{Tikhonov-Fenichel parameter} (TFP); Walcher and collaborators were the first to conceptualize TFPs~\cite{Noethen2011, GOEKE20171PhysicaD,GOEKE1,Goeke2012}.  A manifold of equilibrium points is a special case of a normally hyperbolic manifold; the 
formal definition of a normally hyperbolic invariant manifold can by found in \cite{kuehn2015multiple}. The informal definition of a normally hyperbolic invariant manifold is that the linearized flow that is tangent to the manifold is \textit{dominated} by the flow 
that is normal (transversal) to the manifold. This is more or less intuitive when it comes to manifolds comprised of equilibria, since the flow on such a manifold is trivial. Nevertheless, manifolds of equilibria can lose normal hyperbolicity and, as we will see, the loss can play a significant role in determining when QSS approximations are valid.

The challenging part in the direct utilization of Tikhonov-Fenichel parameters is that often times the perturbed mass action equations 
obtained by regulating a Tikhonov-Fenichel parameter do not necessarily
result in a set of equations that are in Tikhonov standard form.
This can make computing a QSS reduction (without resorting to scaling) 
less than obvious. For example, the parameters $e_0$ and $k_1$ are
Tikhonov-Fenichel parameters that warrant the validity of the 
sQSSA~\cite{GOEKE1}. Treating $k_1$ as a small parameter, we can 
express (\ref{eq:1a})--(\ref{eq:1b}) as
\begin{subequations}
\begin{align}
    \dot{s} &= -\varepsilon k_1^{\star}(e_0-c)s-k_{-1}c\label{eq:ns1},\\
    \dot{c} &= \varepsilon k_1^{\star}(e_0-c)s+(k_{-1}+k_2)c\label{eq:ns2}.
\end{align}
\end{subequations}
The $(s,c)$ phase-plane is equipped with an invariant manifold of 
equilibria, $c=0$, when $k_1=0$. If $0< k_1\ll 1$, then one would 
expect that $s_0\sim s+p$ on the slow timescale, provided the curve 
$c=0$ is normally hyperbolic and attracts nearby trajectories. 
However, determining 
a QSS approximation for $s$ or $c$ from  (\ref{eq:ns1})--(\ref{eq:ns2}) is not obvious without scaling, since the system is not in Tikhonov standard form. 
Nonetheless, using algebraic methods, Walcher and 
collaborators~\cite{Goeke2012,Goeke2014,GoekeBook} developed a 
surprisingly straightforward method for obtaining QSS approximations 
from non-standard systems such as (\ref{eq:ns1})--(\ref{eq:ns2}). In so doing, they have circumnavigated the need for scaling analysis to justify the QSSA. 

For simple planar systems, Tikhonov-Fenichel parameters provide direct insight into the 
conditions that contribute to QSS dynamics, as well the \textit{origin} (i.e., the critical manifold) of 
the corresponding slow manifold, and in this sense they are 
invaluable. However, in order to determine the validity of a 
particular QSS reduction so that experiments can be prepared with 
the intention of determining enzyme activity, the necessary 
``smallness'' of Tikhonov-Fenichel parameters must be defined. 
For example, if $s_0\geq e_0$, then the critical set obtained 
by setting $\varepsilon_1=\varepsilon_2 =0$ in (\ref{eq:sp12})--(\ref{eq:sp22}) yields
\begin{equation}
    e_0s-(s_0-s-p)s =0,
\end{equation}
since the term $e_0-(s_0-s-p)$ can vanish if $e_0 \leq s_0$. In this situation, it is not entirely clear how to go about properly defining 
the critical manifold when $s_0 \geq e_0$. This raises the question 
as to whether or not the rQSSA is still valid when $e_0 \lesssim s_0$ 
and, if so, then how small must the parameters $k_{-1}$ and $k_2$ 
be in order for the approximation~(\ref{eq:rQAAAA}) to hold over 
the slow timescale? Obviously we expect the standard QSSA to be 
valid at extremely large values of $s_0$, but what happens in 
parameter regions where $s_0$ and $e_0$ are the same order of 
magnitude? This question is central not only to the validity of 
the rQSSA, but is also central to understanding the validity of 
the sQSSA. Clearly, it seems that something more than 
scaling analysis and the identification of Tikhonov-Fenichel 
parameters is needed in order to address the validity of the 
rQSSA. We will introduce a method that permits a possible 
solution to the validity of the rQSSA in the subsection that 
follows. 

\subsection{Energy methods: Revisiting the sQSSA and Michaelis--Menten equation} \label{sec:enslavement}


To obtain a fundamental small parameter that ensures the validity of the sQSSA, we need to understand \textit{where} the corresponding slow manifold comes from and, to do this, we must determine the associated critical manifold. One \textit{possible} caveat with the traditional analysis of the sQSSA in the $(s,c)$ phase--plane is that, if we are only interested in bounded trajectories (again, this rules out the particular case of taking $s_0 \to \infty$), then taking $\varepsilon \to 0$ implies that either $k_1$ or $e_0$ must vanish, as Walcher and his collaborators have carefully pointed out \cite{Goeke2012,GOEKE20171PhysicaD,GOEKE1,Noethen2011}. In both of these cases $c$ must vanish, but the dimensionless variable $\bar{c}=c/\varepsilon s_0$ contains a zero denominator in the singular limit. According to Walcher and Lax \cite{Lax2019} this is not a problem. However, in order to avoid having to deal with this, it is useful to write down the dimensionless equations in the $(\bar{s},\bar{p})$ coordinate system:
\begin{subequations}
\begin{align}
\cfrac{\text{d} \bar{s}}{\text{d}\tau} &= -\varepsilon_{SS} \bar{s} + \mu\bar{s}(1-\bar{s}-\bar{p}) + \alpha\beta(1-\bar{s}-\bar{p})\label{eq:Seps}\\
\cfrac{\text{d} \bar{p}}{\text{d}\tau}&=\beta(1-\alpha)(1-\bar{s}-\bar{p})\label{eq:Peps}.
\end{align}
\end{subequations}
Setting $\varepsilon_{SS}=0$ in (\ref{eq:Seps})--(\ref{eq:Peps}) reveals the critical manifold to be $\bar{p}=1-\bar{s}$, which corresponds to the critical manifold of the dimensional system (i.e., $c=0$) obtained by setting $k_1$ or $e_0$ to zero in (\ref{eq:ns1})--(\ref{eq:ns2}). Since the curve $s=s_0-p$ implies $c=0$, we avoid the \textit{other} caveat encountered in the analysis in the $(s,c)$ coordinate system in {{\sc{section \ref{sec:GSPT}}}}: there is no violation of the substrate conservation law $s_0=s+c+p$ along a zeroth-order fast fiber, since initial conditions lie \textit{on} the critical manifold, $c=0$. If a \textit{non-trivial} $c$-nullcline is defined to be the critical manifold, then there is zero loss of substrate along the one-dimensional zeroth-order fibers when $\varepsilon_{SS}=0$, which is a violation of the substrate conservation law, since $s + c > s_0$ as the trajectory approaches a base point that lies above the $s$-axis.

It is expected that trajectories will closely follow the curve $p=s_0-s$ once the perturbation is ``turned on" and $0< \varepsilon_{SS}\ll 1$. To get an idea of how ``good" the approximation $p\sim s_0-s$ is for all time, we need to compute an upper bound ``$M$" such that\footnote{We do not need to write $|\mathcal{E}|$ since $0\leq s_0-s-p$, but later on we will consider cases where the quantity of interest can vary in sign, so it helps to keep the notation consistent through each procedure.}
\begin{equation}
    (1-\bar{s}-\bar{p}) \leq M.
\end{equation}
We will start with a differential equation for the energy\footnote{This is the \textit{mathematical} energy of a function, and should not be confused with the \textit{physical} energy of a particle.} $\mathcal{E}^2 \equiv (s_0-s-p)^2$. Multiplying both sides of the dimensional form of (\ref{eq:Seps}) by ``$(s_0-s-p)\equiv \mathcal{E}$" yields:
\begin{subequations}
\begin{align}
    \cfrac{1}{2}\cfrac{\text{d}\mathcal{E}^2}{\text{d}t} &= k_1e_0s\mathcal{E} - k_1(s+K_M)\mathcal{E}^2,\\
    &\leq k_1e_0s_0|\mathcal{E}| -k_1K_M\mathcal{E}^2.\label{eq:ENERGY}
    \end{align}
\end{subequations}
Now we can can use Cauchy's inequality with ``$\delta$,"
\begin{equation}
    ab \leq \delta a^2 + \cfrac{b^2}{4\delta},
\end{equation}
and expand the term $k_1e_0s_0\mathcal{E}$ in (\ref{eq:ENERGY}):
\begin{equation}
k_1e_0s_0|\mathcal{E}| \leq \cfrac{(k_1e_0s_0)^2}{4\delta} + \delta\mathcal{E}^2.
\end{equation}
Choosing $\delta = k_1K_M/2$ yields
\begin{equation}\label{eq:ENA}
k_1e_0s_0|\mathcal{E}| \leq \cfrac{(k_1e_0s_0)^2}{2k_1K_M} + \cfrac{k_1K_M}{2}\mathcal{E}^2.
\end{equation}
After combining (\ref{eq:ENA}) with (\ref{eq:ENERGY}) we have
\begin{equation}
\cfrac{\text{d}\mathcal{E}^2}{\text{d}t}\leq \cfrac{(k_1e_0s_0)^2}{k_1K_M} -k_1K_M\mathcal{E}^2\label{eq:ENERGYA},
\end{equation}
and integrating both sides of (\ref{eq:ENERGYA}) yields
\begin{equation}\label{eq:supSQ}
    \mathcal{E}^2(t) \leq \mathcal{E}^2(0)e^{\displaystyle-(k_1K_M)t} + s_0^2\eta^2(1-e^{-\displaystyle k_1K_M t}).
\end{equation}
Replacing ``$(1-e^{-k_1K_Mt})$" with ``$1$" and dividing both sides of (\ref{eq:supSQ}) by $s_0^2$ reveals
\begin{equation}\label{eq:LONGTIME}
    \displaystyle \limsup_{t\to\infty} \mathcal{\bar{E}}^2 \leq \eta^2, \quad \mathcal{\bar{E}}\equiv \mathcal{E}/s_0,
\end{equation}
where the limit supremum or ``$\limsup$" of a bounded function $w(t)$ is given by:
\begin{equation}
    \limsup_{t\to \infty} w(t) = \displaystyle \lim_{T\to \infty} \bigg (\sup_{t>T}w(t)\bigg).
\end{equation}
It is clear from (\ref{eq:LONGTIME}) that the long-time error in the approximation is controlled by $\eta$. We can simplify things further: utilizing the fact that $\mathcal{E}(0)=0$, and taking the square root of both sides yields
\begin{equation}\label{eq:limitsup}
(1-\bar{s}-\bar{p}) \leq \eta.
\end{equation}

The bound (\ref{eq:limitsup}) obtained from the energy analysis is promising, but what we are ultimately after is an expression for $p(t)$ that can be utilized for parameter estimation. In practice, the approximation\footnote{This approximation is easily justified through scaling analysis. See \cite{KIM} for details.}
\begin{equation}\label{eq:lambert}
    \cfrac{\text{d}p}{\text{d}t} \sim \cfrac{k_2e_0(s_0-p)}{K_M+(s_0-p)},
\end{equation}
is often used to estimate $K_M$ and $k_2$; the assumption in (\ref{eq:lambert}) is that $c$ is negligible and is well approximated by \begin{equation}\label{eq:cred}
    c\approx \cfrac{e_0(s_0-p)}{K_M+s_0-p}.
\end{equation}
The question that immediately follows is: how \textit{good} is the approximation (\ref{eq:cred})? In order to assess the validity of (\ref{eq:cred}), we will derive an upper bound on $|\mathcal{E}_c|$,
\begin{equation}
\mathcal{E}_c \equiv \bigg(c-\cfrac{e_0(s_0-p)}{K_M+s_0-p}\bigg).
\end{equation}
To begin,
\begin{subequations}
\begin{align}
\cfrac{1}{2}\cfrac{\text{d}\mathcal{E}_c^2}{\text{d}t} &= \bigg(c-\cfrac{e_0(s_0-p)}{K_M+s_0-p}\bigg) \bigg(\cfrac{\text{d}c}{\text{d}t} +\cfrac{e_0K_M}{(K_M+s_0-p)^2}\cfrac{\text{d}p}{\text{d}t}\bigg),\\
&= \mathcal{E}_c \bigg(k_1e_0(s_0-p) -k_1c(e_0+K_M+s_0-p)+k_1c^2 + \cfrac{e_0K_Mk_2c}{(K_M+s_0-p)^2}\bigg),\\
&= \mathcal{E}_c \bigg(-k_1(K_M+s_0-p) \mathcal{E}_c +k_1c(c-e_0) + \cfrac{e_0K_Mk_2c}{(K_M+s_0-p)^2}\bigg),\\
&= -k_1(K_M+s_0-p)\mathcal{E}_c^2 + \mathcal{E}_c\bigg(k_1c(c-e_0) + \cfrac{e_0K_Mk_2c}{(K_M+s_0-p)^2}\bigg)\label{eq:Erhs}.
\end{align}
\end{subequations}
Next, we need to find an upper bound on the right and side of (\ref{eq:Erhs}):
\begin{subequations}
\begin{align}
\cfrac{1}{2}\cfrac{\text{d}\mathcal{E}_c^2}{\text{d}t} &\leq -k_1K_M\mathcal{E}_c^2 + |\mathcal{E}_c|\max \bigg|k_1c(c-e_0) + \cfrac{e_0K_Mk_2c}{(K_M+s_0-p)^2}\bigg|,\\
&\leq -k_1K_M\mathcal{E}_c^2  + |\mathcal{E}_c|\max|k_1c(c-e_0)|+|\mathcal{E}_c|\max \bigg|\cfrac{e_0K_Mk_2c}{(K_M+s_0-p)^2}\bigg|,\\
&\leq -k_1K_M\mathcal{E}_c^2 + |\mathcal{E}_c|\bigg(\cfrac{k_1e_0^2}{4} + \cfrac{e_0k_2\lambda}{K_M}\bigg).
\end{align}
\end{subequations}
Using Cauchy's inequality with $\delta = k_1K_M/2$, it follows that
\begin{equation}\label{eq:cenergy}
  \cfrac{\text{d}\mathcal{E}_c^2}{\text{d}t} \leq -k_1K_M \mathcal{E}_c^2 +  \cfrac{1}{k_1K_M}\bigg(\cfrac{k_1e_0^2}{4} + \cfrac{e_0k_2\lambda}{K_M}\bigg)^2. 
\end{equation}
Applying Gronwall's lemma and integrating both sides of (\ref{eq:cenergy}) yields
\begin{equation}\label{eq:expo}
\mathcal{E}_c^2 \leq \mathcal{E}_c^2(0)\displaystyle e^{\displaystyle -k_1K_Mt} + \cfrac{1}{k_1^2K_M^2}\bigg(\cfrac{k_1e_0^2}{4} + \cfrac{e_0k_2\lambda}{K_M}\bigg)^2.
\end{equation}
Although the exponential term (\ref{eq:expo}) is in terms of the dimensional time, $t$, we expect that if  (\ref{eq:expo}) is written in terms of the slow time, $\bar{T} = tk_2\varepsilon_{SS}$, then it should decay relatively quickly in $\bar{T}$. To do this we need to solve
\begin{equation}\label{eq:timeconstant}
    k_1K_Mt = C_{\bar{T}} k_2\varepsilon_{SS}t
\end{equation}
for the constant, $C_{\bar{T}}$. Solving (\ref{eq:timeconstant}) reveals $C_{\bar{T}}=1/\nu \varepsilon$, and it clear that the exponential term in (\ref{eq:expo}) will vanish rapidly in $\bar{T}$ whenever $\varepsilon\nu \ll 1$. Finally, taking the square root of both sides while noting that $\sqrt{x+y}\leq \sqrt{x}+\sqrt{y}\;\;\forall x, y \geq 0$, and dividing through by $e_0$ produces the following upper bound:
\begin{equation}\label{eq:limitsupA}
 |\tilde{c} -h(\bar{p})|\leq \mu \displaystyle e^{\displaystyle -\bar{T}/2\varepsilon_{SS}\nu}   +\cfrac{\eta}{4} + \nu\cfrac{\lambda}{K_M},  \quad h(\bar{p}) \equiv \cfrac{\sigma(1-\bar{p})}{1+\sigma(1-\bar{p})},
\end{equation}
where $\tilde{c}\equiv c/e_0$. The bound (\ref{eq:limitsupA}) is actually quite revealing. First, we immediately see that taking $\nu \to 0$ does not support the validity of the sQSSA (i.e., the MM equation), even though the exponential term decays rapidly in the slow time whenever $\varepsilon_{SS}\nu \ll 1$, and ``nullcline following" in $(s,c)$ phase-plane is observed when $k_2 \ll k_{-1}$. Thus, we do not encounter any sort of dilemma between the sQSSA and the extended QSSA when energy methods are appropriately employed. Second, if $s_0 < e_0$, then $\lambda < s_0$, and we have
\begin{equation}
 |\tilde{c} -h(\bar{p})|\leq \mu \displaystyle e^{\displaystyle -\bar{T}/2\varepsilon_{SS}\nu}   +\cfrac{\eta}{4} + \nu\sigma.  
\end{equation}
Since $\sigma < \eta$ when $s_0 < e_0$, it follows that $\eta \to 0$ in order for the sQSSA to be valid. Third, note that making $s_0$ arbitrarily large only makes the exponential term decay faster, since $\varepsilon_{SS}\to0$ as $s_0 \to \infty$. However, large $s_0$ has no influence on the long-time bound: only $K_M$ and $e_0$ influence the long-time bound. Finally, if $e_0 < s_0$, then $\lambda \leq e_0$ and it holds that
\begin{equation}
  |\tilde{c} -h(\bar{p})|\leq \mu \displaystyle e^{\displaystyle -\bar{T}/2\varepsilon_{SS}\nu}   +\cfrac{5}{4}\eta. \phantom{+} \phantom{\nu\sigma} 
\end{equation}

\begin{figure}[hbt!]
\centering
\includegraphics[width=10cm]{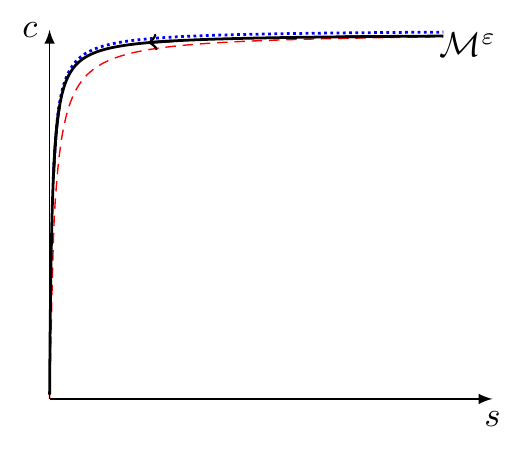}
\caption{\textbf{Near the origin, the slow manifold lies close to the $s$-nullcline in the $(s,c)$ phase--plane when $1\ll \eta$}.
This figure is an illustration and is meant to convey conceptual features that occur in the phase--plane when $\eta \gg 1$. The solid black curve is the slow manifold; the dotted blue curve is the $s$-nullcline, and the dashed red curve is the $c$-nullcline. When $\eta$ is sufficiently large, the slow manifold lies very close to the $s$-nullcline at sufficiently small $s$. Comparing with {{\sc{figure \ref{fig:11}}}}, this is exactly the condition we need in order to apply the rQSSA}\label{fig:sNULL}
\end{figure}

In all cases the natural ``small parameter" that arises from the energy analysis is $\eta$. What can happen if $\varepsilon_{SS}\ll 1$ but $\eta \gg 1$? The answer is that a curious situation arises: the $(s,c)$ phase-plane trajectory will rapidly approach the $c$-nullcline, and proceed to follow it closely for a finite amount of time. However, as the reaction progresses, the phase-plane trajectory, which follows the invariant manifold $\mathcal{M}^{\varepsilon}$, will start to move away from the $c$-nullcline, and eventually follow the $s$-nullcline (see {\sc figure \ref{fig:sNULL}}).

In conclusion of this section, our analysis suggests that $\eta$ is 
the fundamental small parameter that justifies the validity of the 
sQSSA. The condition $\eta \ll 1$ can be viewed as a special 
case of $\varepsilon_{SS} \ll 1$. However, if $\varepsilon_{SS} \ll 1$ 
and $1 \ll \eta$, then we observe a scenario in the $(s,c)$ phase-plane 
in which the sQSSA is initially valid after the fast transient, 
but steadily loses validity until the rQSSA becomes valid for the 
remainder of the reaction. It seems as though the situation is quite unique, since the sQSSA appears to valid over the slow timescale, 
but the rQSSA is valid over a super-slow timescale. If we can understand 
the origin of the slow manifold in this case (i.e., its associated 
critical manifold), then perhaps it will be possible to explain why this 
scenario occurs. Of course, this raises the question: ``Can a unique critical 
manifold be determined in this case?'' The answer to this question resides in the 
analysis of the rQSSA, and we will analyze this particular approximation in great detail in {{\sc section \ref{sec:rQSSA}}}. 

\section{The reverse Quasi-Steady-State Approximation} \label{sec:rQSSA}
In this section, we derive the validity of the rQSSA using a combination 
of scaling and energy methods.


\subsection{Justification of the reverse Quasi-Steady-State Approximation
in the $(p,c)$ phase-plane: Scaling and bifurcation analysis}\label{sec:huer}
\begin{figure}[hbt!]
\centering
\includegraphics[width=8cm]{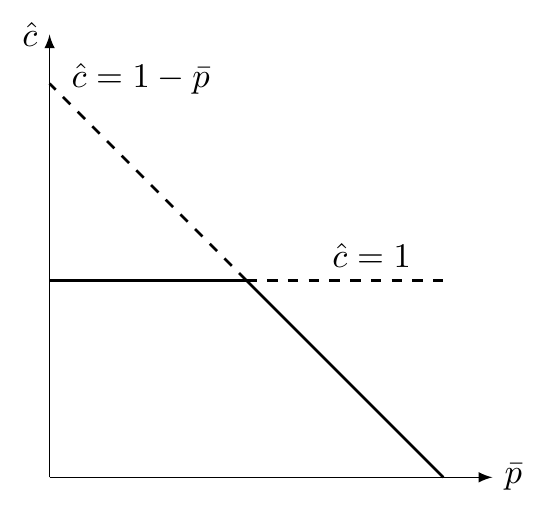}
\includegraphics[width=8cm]{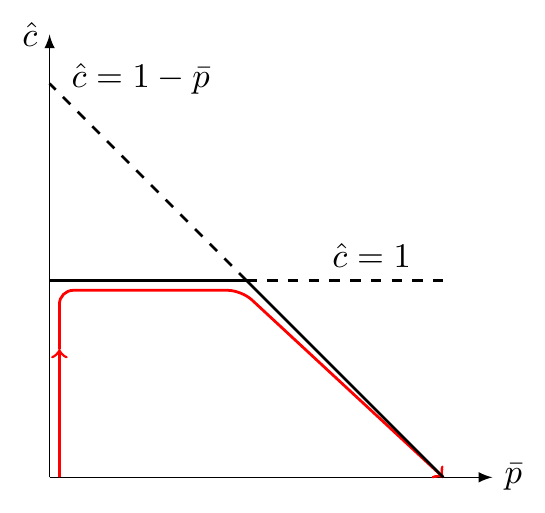}
\caption{\textbf{The critical set of the rQSSA when $e_0<s_0$ 
and $K_M=0$}. This figure provides a visualization of 
the invariant set recovered by setting $e_0<s_0$ and $k_{-1}=k_2=0$. {\sc Left}: The critical set contains two 
transversely intersecting branches of fixed points. Thick solid 
curves correspond to non-isolated \textit{stable} fixed points, and 
thick dashed lines correspond to \textit{unstable} fixed points. 
The horizontal curve corresponds to the critical set $\hat{c}=1$, 
and the diagonal curve corresponds to the critical set 
$\bar{p}+\hat{c}=1$. The trajectory rapidly approaches the curve 
$\hat{c}=1$, then reaches the curve $1=\hat{c}+\bar{p}$ once 
$t\sim t_{\ell}$, and begins descending towards the origin. Clearly, 
each set constitutes a normally hyperbolic manifold 
everywhere except where the branches intersect, which 
corresponds to a transcritical singularity. {\sc Right}: A 
typical trajectory (red solid curve) closely follows the 
attracting critical submanifolds once the perturbation is turned 
on.}
\label{fig:A}
\end{figure}

It is straightforward to verify the validity of the rQSSA through scaling analysis. 
Assuming that $s_0 < e_0$, it is obvious that $c$ is bounded by 
$s_0$, which suggests that $\hat{c} = c/s_0$ is the appropriate scaled 
variable when $s_0 < e_0$. Rescaling the mass action equations (\ref{eq:tQSSA10})--(\ref{eq:tQSSA20}) in the
$(p,c)$ phase--plane yields
\begin{subequations}
\begin{align}
 \varepsilon^*\nu\cfrac{\text{d} \hat{c}}{\text{d}\tilde{T}} &= (1-\ell\hat{c})(1-\hat{c}-\bar{p}) - \varepsilon^*\hat{c},\\
  \cfrac{\text{d} \bar{p}}{\text{d}\tilde{T}} &= \hat{c},
\end{align}
\end{subequations}
and the corresponding critical manifold, $\mathcal{M}_0$, is obtained by setting $\varepsilon^*=0$,
\begin{equation}
\mathcal{M}_0 = \{(\bar{p},\hat{c})\in \mathbb{R}^2_+ | (1-\ell \hat{c})(1-\hat{c}-\bar{p})=0 \}.
\end{equation}
If $s_0 < e_0$, then $(1-\ell \hat{c}) >0$, and we can take the 
critical manifold to be $1-\hat{c}-\bar{p}=0$. Furthermore, the critical manifold is normally hyperbolic and attracting, since
\begin{equation}
    \cfrac{\partial}{\partial \hat{c}}(1-\ell \hat{c})( 1-\hat{c}-\bar{p}) <-.
\end{equation}
Thus, due to the fact that the critical manifold is both normally hyperbolic and attracting, we expect trajectories to approach the curve $c=s_0-p$ when $s_0< e_0$ and $0 < \varepsilon^* \ll 1$. But what happens when $e_0\leq s_0$? Let us start by setting $e_0=s_0$, in which case $\ell =1$. 
The dimensionless equations are given by
\begin{subequations}
\begin{align}
 \varepsilon^*\nu\cfrac{\text{d} \hat{c}}{\text{d}\tilde{T}} &= (1-\hat{c})(1-\hat{c}-\bar{p}) - \varepsilon^*\hat{c},\\
  \cfrac{\text{d} \bar{p}}{\text{d}\tilde{T}} &= \hat{c},
\end{align}
\end{subequations}
and the associated critical \textit{set} is
\begin{equation}\label{eq:SET}
\mathcal{M}_0 = \{(\bar{p},\hat{c})\in \mathbb{R}^2_+| (1-\hat{c})(1-\hat{c}-\bar{p})=0\}.
\end{equation}
The set (\ref{eq:SET}) fails to be a manifold at the 
point where the curves $\hat{c}=1$ and $\hat{c} = 1-\bar{p}$ intersect, 
which is precisely the point $(\bar{p},\hat{c})=(0,1)$. Moreover, 
there is a loss of normal hyperbolicity at this point, since
\begin{equation}
\cfrac{\partial}{\partial \hat{c}}\bigg[(1-\hat{c})(1-\hat{c}-\bar{p})\bigg]\bigg|_{(\bar{p},\hat{c})=(0,1)}=0.
\end{equation}
What is happening at the point where normal hyperbolicity is lost? The sub-manifolds $\hat{c}=1$ and $\bar{p}=1-\hat{c}$ intersect and exchange stability at the point $(\bar{p},\hat{c})=(0,1)$ when $K_M=0$. To see this, rescale time and define $\tau^* = \tilde{T}/\varepsilon^*\nu$:
\begin{subequations}
\begin{align}
 \cfrac{\text{d} \hat{c}}{\phantom{i}\text{d}\tau^*} &= (1-\hat{c})(1-\hat{c}-\bar{p}) - \varepsilon^*\hat{c},\\
  \cfrac{\text{d} \bar{p}}{\phantom{i}\text{d}\tau^*} &= \varepsilon^*\nu\hat{c}.
\end{align}
\end{subequations}
Since $e_0 = s_0$, setting $\varepsilon^*=0$ implies that $K_M=0$, 
and we recover
\begin{subequations}
\begin{align}
 \cfrac{\text{d} \hat{c}}{\phantom{i}\text{d}\tau^*} &= (1-\hat{c})(1-\hat{c}-\bar{p}),\label{eq:fastsub1}\\
  \cfrac{\text{d} \bar{p}}{\phantom{i}\text{d}\tau^*} &= 0\label{eq:fastsub2}.
\end{align}
\end{subequations}
Treating $\bar{p}$ as a slowly-varying parameter 
in the fast subsystem given by (\ref{eq:fastsub1})--(\ref{eq:fastsub2}), 
the normal form of the transcritical bifurcation 
is recovered by making the change of variables $u = 1-\hat{c}$
\begin{equation}\label{eq:NF}
\frac{\text{d}u}{\phantom{i}\text{d}\tau^*} = \bar{p}u - u^2.
\end{equation}
The normal form equation (\ref{eq:NF}) indicates that the critical submanifolds\footnote{The transcritical bifurcation is also easily recoverable in $(s,p)$ coordinates.} undergo an exchange of stability at the transcritical singularity. Clearly, this transcritical singularity influences the dynamics whenever $K_M=0$ and $e_0\leq s_0$ ({see {\sc Figure \ref{fig:A})}}. 

\subsection{A heuristic analysis based on timescale separation} 
\begin{figure}[hbt!]
\centering
\includegraphics[width=8cm]{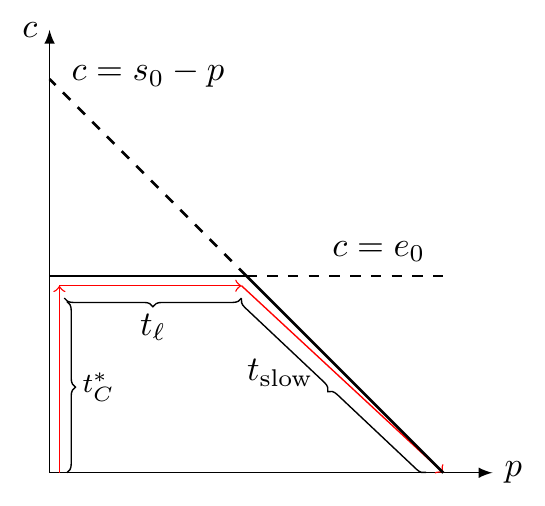}
\caption{\textbf{There are three stages in the dynamics of the MM reaction mechanism (\ref{eq:mm1}) when $e_0 < s_0$ and 
$K_M \approx 0$}. The red arrows demarcate the approximate path of a typical trajectory during each stage. Stage one is the fast transient: the trajectory rapidly approaches the curve $c=e_0$, and the duration of this stage is roughly $t_C^*$. During the second stage, the trajectory closely follows the curve $c=e_0$; the duration of this stage is roughly $t_{\ell}$. The final stage corresponds to the the rQSSA, as the trajectory follows the curve $c=s_0-p$. The timescale corresponding to stage three is $t_{\text{slow}}=1/k_2$. If $t_C^* \ll t_{\ell} \ll 1/k_2$, then the rQSSA is at least heuristically valid.}
\label{fig:huer}
\end{figure}

The bifurcation analysis from the previous section gave us a global picture of what the phase-plane dynamics look like when $K_M=0$. Moreover, we can use the bifurcation structure to our advantage and attempt to determine a suitable ``small parameter" that warrants the validity of the rQSSA. To begin, if $e_0 < s_0$, then there are essentially three stages to the reaction (see {{\sc figure \ref{fig:huer}}}).

The first stage is the transient stage, and the phase-plane trajectory moves almost vertically towards the curve $c=e_0$. The approximate duration of this stage is given by Tzafriri's fast timescale, $t_C^*$. After the initial fast transient the differential equation for product 
formation is roughly
\begin{equation}
\dot{p} \sim k_2e_0.
\end{equation}
Integrating the above expression yields $p\sim k_2e_0t$. 
Since $c\approx e_0$ in this regime, the trajectory will reach the 
vicinity of the $c$-axis once $p=s_0-e_0$, which, based on the
substrate conservation law, would imply that $s\approx 0$. Thus, 
the time it takes the trajectory to reach the curve $s_0=c+p$ is given 
by $t_{\ell}$:
\begin{equation}
t_{\ell} = \cfrac{s_0-e_0}{k_2e_0}.
\end{equation}
Next, the trajectory can be assumed 
to be very close to the curve $c=s_0-p$ and, at this point, the time it 
takes for the reaction to effectively complete is $t_{\text{slow}}=1/k_2$. If 
the ratio of $t_{\ell}$ to $t_{\text{slow}}$ is small, i.e.,
\begin{equation}\label{eq:frac}
0<\cfrac{t_{\ell}}{t_{\text{slow}}} =\cfrac{s_0-e_0}{e_0}\ll 1,
\end{equation}
then, from a \textit{practical} point of view (as opposed to a rigorous point of view), the rQSSA can be 
considered valid if $s_0>e_0$ and $K_M$ is incredibly small. 

The heuristic analysis at least partially supports the idea that parameter domain over which the rQSSA is valid may in fact be larger than previously assumed. However, the combined conditions, $\varepsilon^* \ll 1$ and $(s_0-e_0)/e_0 \ll 1$, still do not explain how 
increasing $s_0$ impacts the validity of the approximation. This 
is a reasonable question that has plagued theorists for some time now (see, for example~\cite{Segel1989} and~\cite{SCHNELL2000}). The usual methods of scaling analysis and timescale estimation have not provided very satisfying answers. Consequently, this suggests that a new method needs to be employed to tackle this problem. 

\subsection{The reverse quasi-steady-state approximation: energy methods}
In order to determine a small parameter that ensures the validity of the rQSSA, it helps to first consider what we expect to happen \textit{chemically} when such a parameter is made arbitrarily small. Historically, thanks to the careful work of Nguyen and Fraser \cite{HNsF}, the geometric interpretation of the rQSSA's validity has been associated with a phase-plane trajectory closely following the $s$-nullcline in the $(s,c)$ phase--plane. However, instead of thinking about the validity of the rQSSA in terms of a trajectory's propensity for following the $s$-nullcline, it is perhaps better to think of it in purely chemical terms: when the rQSSA is valid, the initial substrate concentration vanishes rapidly with respect to the timescale $\tilde{T}\equiv k_2t$. Translated mathematically, this means that $s$ must \text{dissipate} rapidly with respect to $\tilde{T}$. Thus, we want an equation that describes the dissipation of $s^2/2$:
\begin{subequations}
\begin{align}
    \cfrac{1}{2}\cfrac{\text{d} s^2}{\text{d}t}&= -k_1(e_0-c)s^2 +k_{-1}cs,\\
    &\leq -k_1(e_0-\lambda)s^2 + k_{-1}s\lambda.
\end{align}
\end{subequations}
The term $k_{-1}s\lambda$ can be expanded thanks to Cauchy's inequality
\begin{equation}
    k_{-1}s\lambda \leq \delta s^2 + \cfrac{k_{-1}^2\lambda^2}{4\delta},
\end{equation}
and choosing $\delta = k_1(e_0-\lambda)/2$  yields
\begin{equation}
\cfrac{1}{2}\cfrac{\text{d} s^2}{\text{d}t}\leq -\cfrac{1}{2}k_1(e_0-c)s^2 +\cfrac{k_{-1}^2\lambda^2}{2k_1(e_0-\lambda)}.
\end{equation}
From Gronwall's lemma we have\footnote{From this point forward, we will drop exponential terms of the form ``$(1-e^{-\lambda t})$" that converge to $1$ as $t\to \infty$ and just replace it by ``$1$."}
\begin{subequations}
\begin{align}
    s^2(t) &\leq s_0^2 e^{\displaystyle -k_1(e_0-\lambda)t} + \cfrac{k_{-1}^2\lambda^2}{k_1^2(e_0-\lambda)^2}\bigg(1-e^{\displaystyle -k_1(e_0-\lambda)t}\bigg),\\
    &\leq s_0^2 e^{\displaystyle -k_1(e_0-\lambda)t} + \cfrac{k_{-1}^2\lambda^2}{k_1^2(e_0-\lambda)^2}.
    \end{align}
\end{subequations}
Diving both sides by $s_0^2$, taking the square root of both sides, and invoking the triangle inequality yields:
\begin{equation}\label{eq:senergy}
    s(t) \leq s_0e^{\displaystyle -k_1(e_0-\lambda)t/2} + \cfrac{k_{-1}\lambda}{k_1(e_0-\lambda)}.
\end{equation}
Next, re-writing (\ref{eq:senergy}) in terms of $\tilde{T}$,
\begin{subequations}\label{eq:dsenergy}
\begin{align}
 \bar{s}(\tilde{T}) &\leq e^{\displaystyle -\tilde{T}/2\underline{\varepsilon}\nu} + \tilde{\varepsilon}\cfrac{\lambda}{s_0},\\
 &\leq e^{\displaystyle -\tilde{T}/2\underline{\varepsilon}\nu} + \tilde{\varepsilon},
 \end{align}
\end{subequations}
reveals two small parameters:
\begin{equation}
    \underline{\varepsilon} \equiv \cfrac{K_M}{e_0-\lambda}, \qquad \tilde{\varepsilon}\equiv \cfrac{K_S}{e_0-\lambda}.
\end{equation}
Now we have a choice to make: which small parameter corresponds to the potential validity of the rQSSA? The answer to this question resides in Fenichel theory: Setting $\underline{\varepsilon}$ to zero results in a critical manifold of fixed points, whereas making $\tilde{\varepsilon}=0$ by setting $k_{-1}=0$ results in the manifold $s_0=c+p$ being invariant, but it is not a manifold of equilibrium points, since $k_{-1}$ is not a TFP (see \cite{GOEKE20171PhysicaD} for details). Moreover, it holds that
  \begin{equation}
   \bar{s}(\tilde{T}) \leq e^{\displaystyle -\tilde{T}/2\underline{\varepsilon}} + \underline{\varepsilon},
\end{equation}
since $\nu \leq 1$ and $ \tilde{\varepsilon} \leq \underline{\varepsilon}$. Consequently, we expect $s$ to vanish rapidly in $\tilde{T}$ whenever $\underline{\varepsilon}\ll 1.$ Furthermore, note that
\begin{equation}
    \lim_{K_M \to 0} \underline{\varepsilon} = 0 \;\;\text{if}\;\; s_0 < e_0, \qquad \lim_{K_M \to 0} \underline{\varepsilon} = \cfrac{s_0-e_0}{e_0}\;\;\text{if}\;\;s_0 > e_0,
\end{equation}
and we recover the exact small parameter (\ref{eq:frac}) that was obtained heuristically in Section \ref{sec:huer}. 

As a final check, we observe the influence of $\underline{\varepsilon}$ on the dynamics in the $(p,c)$ phase-plane when $0 < \underline{\varepsilon}\ll 1$. First note that the $c$-nullcline can be approximated by the curve $s_0=c+p$ as $K_M \to 0$ in any region where $s_0 < e_0$. Thus, if trajectories closely follow the $c$-nullcline when $K_M \ll 1$ and $s_0 < e_0$, then the rQSSA should be valid. Second, one can show that
\begin{equation}\label{eq:finalbounder}
\displaystyle |c-c_0(p)| \leq \lambda e^{\displaystyle-t\zeta_T/2} + \cfrac{\eta\lambda}{(1+\eta)(1+\kappa)}\bigg(\cfrac{K_M}{e_0-\lambda + K_M}\bigg),
\end{equation}
where $\zeta_T \equiv k_1(e_0-\lambda + K_M)$ (we will derive this in the Section \ref{sec:timescales} using energy methods, but also in Section \ref{sec:appen} using a method pioneered by Gradshtein \cite{Gradshtein1953}). Dividing both sides by $\lambda$ and expressing (\ref{eq:finalbounder}) in terms of $\tilde{T}$ yields
\begin{equation}\label{eq:finalbounder1}
\displaystyle |\bar{c}-\bar{c}_0(\bar{p})| \leq e^{\displaystyle-\tilde{T}/2\underline{\varepsilon}\nu} + \cfrac{\eta\nu\underline{\varepsilon}}{(1+\eta)(1+\underline{\varepsilon})},
\end{equation}
from which it follows that
\begin{equation}\label{eq:finalbounder11}
\displaystyle |\bar{c}-\bar{c}_0(\bar{p})| \leq  e^{\displaystyle-\tilde{T}/2\underline{\varepsilon}}+ \underline{\varepsilon}.
\end{equation}
The upper bound (\ref{eq:finalbounder11}) reveals that the phase-plane trajectory will closely follow the $c$-nullcline after a brief transient. Since the $c$-nullcline is approximately $s_0 =c+p$ whenever $s_0<e_0$ and $\underline{\varepsilon}\ll 1$, it holds that the rQSSA is valid under these conditions, and we take $\underline{\varepsilon}$ to be \textit{the} small parameter that ensures its validity. 


In summary of this section, we have shown that, while scaling 
methods or Tikhonov-Fenichel parameters can be used to predict 
the presence of a singular point, more analysis is needed to assess the validity of the rQSSA in neighborhoods impacted 
by the critical singularity. Hence, the presence of the singular point in the critical set is what makes the analysis of the rQSSA  
challenging. As we have shown, energy methods can be easily employed to assess the validity of 
the rQSSA when initial substrate and enzyme concentrations are of 
similar magnitude. Moreover, the bound 
\begin{equation}\label{eq:finalbound}
\displaystyle |c-c_0(p)| \leq \lambda e^{\displaystyle-t\zeta_T/2} + \cfrac{\eta\lambda}{(1+\eta)(1+\kappa)}\bigg(\cfrac{K_M}{e_0-\lambda + K_M}\bigg),
\end{equation}
provides direct insight into the origin of each critical manifold, and suggests that there 
are three parameters, $\eta$, $\underline{\varepsilon}$ and $\nu$, 
that, when adequately small, guarantee the long-time validity of the QSSA, the 
rQSSA and the tQSSA, respectively. We will discuss this observation in more detail in {\sc{section \ref{sec:timescales}}}.

\section{Final remarks on timescales and small parameters}\label{sec:timescales}

In this section we explore the relationship between the fast and slow timescales of the tQSSA, Tikhonov-Fenichel parameters, and the small parameters obtained from the energy analysis. 

\subsection{Timescale separation and small parameters: local versus global conditions for the accuracy of QSS approximations}

As mentioned earlier, because the tQSSA encompasses the rQSSA and the sQSSA, the accepted criterion for the validity of a QSSA is separation of Tzafriri's timescales, $\varepsilon_{T}\ll 1$ (see \ref{eq:TZA}). We can express the bound on the $\limsup$ of $\underline{c}$ and $\overline{c}_0(\bar{p})$ in terms of $T_z = t/t_P$,
\begin{equation}\label{eq:totalLIMSUP}
|\underline{c}-\bar{c}_0(\bar{p})| \leq e^{\displaystyle-T_z/2\varepsilon_D} + \varepsilon_{L}.
\end{equation}
To compute (\ref{eq:totalLIMSUP}), we start by computing,
\begin{equation}
\displaystyle \limsup_{t \to \infty} (c-h^-(p))^2,
\end{equation}
where we will use the notation
\begin{subequations}
\begin{align}
h^-(p) &= \cfrac{1}{2}(e_0+s_0+K_M-p)-\cfrac{1}{2}\sqrt{(e_0+s_0+K_M-p)^2-4e_0(s_0-p)},\\
h^+(p) &= \cfrac{1}{2}(e_0+s_0+K_M-p)+\cfrac{1}{2}\sqrt{(e_0+s_0+K_M-p)^2-4e_0(s_0-p)},
\end{align}
\end{subequations}
for convenience. Thus, the reader should identify $\lambda = h^-(0)$ and note that $c_0(p)\equiv h^-(p)$. The differential equation for the energy, $\mathcal{E}^2\equiv (c-h^-(p))^2$, is
\begin{equation}
\cfrac{1}{2}\cfrac{\text{d}\mathcal{E}^2}{\text{d}t} = \bigg[\cfrac{\text{d}c}{\text{d}t}-\cfrac{\text{d}h^-}{\text{d}p}\cfrac{\text{d}p}{\text{d}t}\bigg]\mathcal{E}.
\end{equation}
Carefully note that the derivative for $c$ can be factored, $\dot{c} = k_1(c-h^-(p))(c-h^+(p))$, and thus
\begin{subequations}
\begin{align}
\cfrac{1}{2}\cfrac{\text{d}\mathcal{E}^2}{\text{d}t} &= \bigg[k_1(c-h^-(p))(c-h^+(p))-\cfrac{\text{d}h^-}{\text{d}p}\cfrac{\text{d}p}{\text{d}t}\bigg]\mathcal{E}\\
&=\bigg[k_1\mathcal{E}(c-h^+(p))-\cfrac{\text{d}h^-}{\text{d}p}\cfrac{\text{d}p}{\text{d}t}\bigg]\mathcal{E}\\
&\leq \mathcal{E}^2k_1(c-h^+(p))+\max\bigg|\cfrac{\text{d}h^-}{\text{d}p}\bigg|\sup\bigg|\cfrac{\text{d}p}{\text{d}t}\bigg||\mathcal{E}|\label{eq:118c}.
\end{align}
\end{subequations}
Next, the term ``$c-h^+(p)$" can be bounded above:
\begin{equation}\label{eq:UP}
    \max c-h^+(p) \leq \lambda -\min h^+(p) = -k_1(e_0+K_M-\lambda).
\end{equation}
Using the Cauchy-$\delta$ inequality with $\delta=k_1(e_0+K_M-\lambda)/2$, we obtain
\begin{equation}\label{eq:BIGeq}
\max\bigg|\cfrac{\text{d}h^-}{\text{d}p}\bigg|\sup\bigg|\cfrac{\text{d}p}{\text{d}t}\bigg||\mathcal{E}| \leq \cfrac{k_1(e_0+K_M-\lambda)}{2}\mathcal{E}^2 + \cfrac{\bigg(\max\bigg|\cfrac{\text{d}h^-}{\text{d}p}\bigg|\bigg)^2\bigg(\sup\bigg|\cfrac{\text{d}p}{\text{d}t}\bigg|\bigg)^2}{2k_1(e_0+K_M-\lambda)}.
\end{equation}
Inserting (\ref{eq:BIGeq}) into (\ref{eq:118c}) yields
\begin{equation}
 \cfrac{\text{d}\mathcal{E}^2}{\text{d}t}   \leq -\mathcal{E}^2k_1(e_0+K_M-\lambda)+\cfrac{\bigg(\max\bigg|\cfrac{\text{d}h^-}{\text{d}p}\bigg|\bigg)^2\bigg(\sup\bigg|\cfrac{\text{d}p}{\text{d}t}\bigg|\bigg)^2}{k_1(e_0+K_M-\lambda)},
\end{equation}
and from Gronwall's Lemma and the triangle inequality we recover:
\begin{equation}\label{eq:totalLIMSUPP}
|\mathcal{E}(t)| \leq |\mathcal{E}(0)|e^{-\displaystyle k_1(e_0+K_M-\lambda)t} + \cfrac{\max\bigg|\cfrac{\text{d}h^-}{\text{d}p}\bigg|\sup\bigg|\cfrac{\text{d}p}{\text{d}t}\bigg|}{k_1(e_0+K_M-\lambda)}.
\end{equation}
Calculating the ``$\max$" and ``$\sup$" on the right hand side of (\ref{eq:totalLIMSUPP}), dividing both sides by $\lambda$, and expressing the exponential in terms of $T_z = t/t_P$ yields (\ref{eq:totalLIMSUP}), where $\varepsilon_D$ denotes the ``exponential $\varepsilon$,'' and $\varepsilon_L$ denotes the ``long-time $\varepsilon$:''
\begin{equation}\label{eq:hbound}
\varepsilon_{D} \equiv \bigg(\cfrac{\lambda}{s_0}\bigg)\bigg(\cfrac{1}{1+\kappa}\bigg)\bigg(\cfrac{\underline{\varepsilon}}{1+\underline{\varepsilon}}\bigg), \qquad \varepsilon_L \equiv \bigg(\cfrac{\eta}{1+\eta}\bigg)\bigg(\cfrac{1}{1+\kappa}\bigg)\bigg(\cfrac{\underline{\varepsilon}}{1+\underline{\varepsilon}}\bigg).
\end{equation}
As expected, $\varepsilon_D$ is dependent on initial data, meaning its magnitude is influenced by the starting location of the trajectory within the vector field. Moreover, $\varepsilon_{D}$ can be ``factored" into three small parameters:
\begin{enumerate}[label=(\roman*)]
    \item $\lambda/s_0$: vanishes as a result of zero enzyme to substrate ratio.\\
    \item $1/(1+\kappa) \equiv \nu$: vanishes in the limit of infinitely slow product formation.\\
    \item $\underline{\varepsilon}/(1+\underline{\varepsilon})$: vanishes at infinitely high enzyme concentration or, zero $K_M$ whenever $s_0\leq e_0$.
\end{enumerate}

\begin{figure}[hbt!]
\centering
\includegraphics[width=11cm]{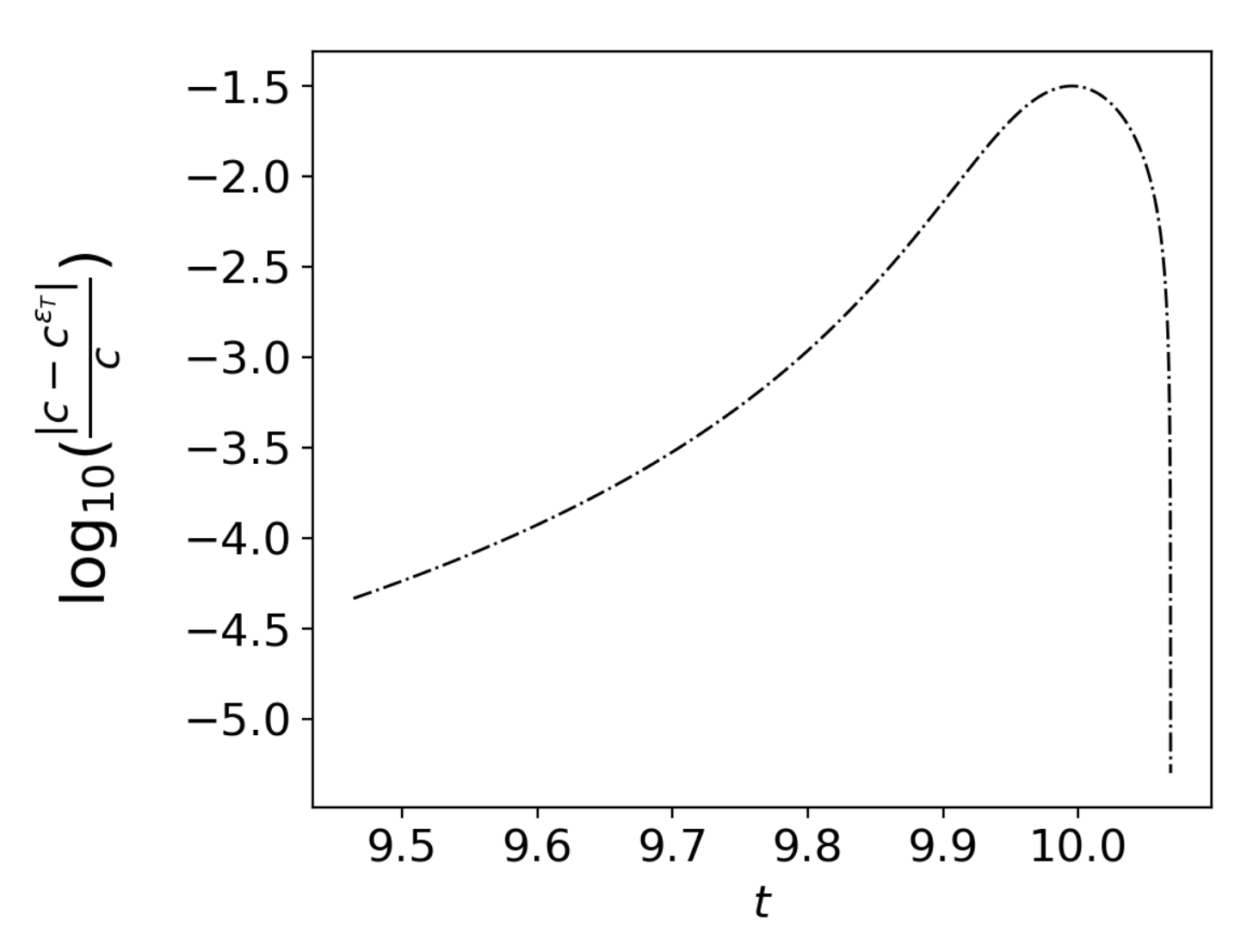}
\includegraphics[width=11cm]{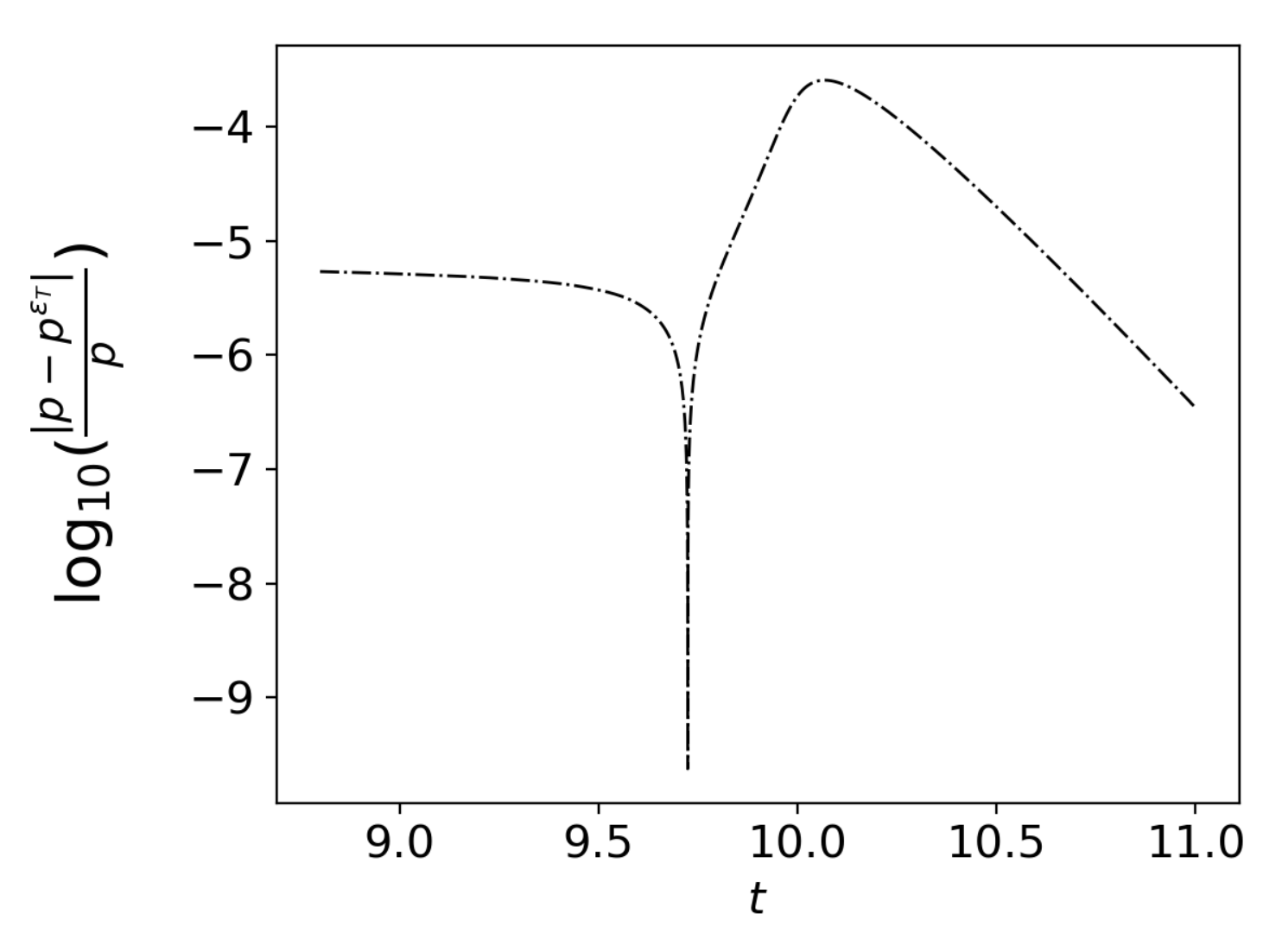}
\caption{\textbf{The tQSSA can worsen as $t\to \infty$ if $\varepsilon_{LT}$ is not sufficiently small}. In both panels, the numerical solution is obtained with parameter values: $e_0=10, s_0 = 1000, k_1=20,k_{-1}=10,k_2=10$. Thus, $\varepsilon_T \sim 10^{-6}$, but $\varepsilon_{LT}\approx 0.1$. {\sc{top panel}}: The $\log_{10}$ of the relative error between the numerical solution to the mass action equations (\ref{eq:tQSSA10})-(\ref{eq:tQSSA20}) for $c$, and the numerical solution to the asymptotic approximation (\ref{eq:mroot}), denotes $c^{\varepsilon_T}$, is demarcated by the dashed-dotted line. {\sc{bottom Panel}}: The dashed-dotted line demarcates the $\log_{10}$ of the relative error between the numerical solution to the mass action equations (\ref{eq:tQSSA10})-(\ref{eq:tQSSA20}) for $p$, and the numerical solution to the asymptotic approximation (\ref{eq:tQ}), denoted $p^{\varepsilon_T}$. Even though timescale separation is quite large, there are regions in the numerically-obtained time course data where $c^{\varepsilon_T}$ approximates $c$ only to within one digit of accuracy. While the loss in accuracy may not be outright disastrous in terms of parameter estimation, it is an observation that is worth pointing out in order to gain a deeper quantitative understanding of the various QSS approximations utilized for Michaelis--Menten reaction mechanism (\ref{eq:mm1}). The parameters have been assigned arbitrary units for illustrative purposes.}
\label{fig:final}
\end{figure}

What is the relationship between Tzafriri's timescales and $\varepsilon_D$ and $\varepsilon_L$? Tzafriri's fast timescale, $t_C^*$, can be written as
\begin{equation}\label{eq:tCiden}
t_C^* = \cfrac{1}{k_1(K_M+e_0-\lambda + s_0-\lambda)}.
\end{equation}
Using the identity (\ref{eq:tCiden}), a simple calculation reveals $\varepsilon_T \leq \varepsilon_D$. In fact, $\varepsilon_T \sim \varepsilon_D$ at very high enzyme concentration, and $\varepsilon_T\to 0$ if $\varepsilon_D \to 0.$ Thus, we can control timescale separation by controlling $\varepsilon_D$. However, $\varepsilon_D \ll 1$ only implies that the phase-plane trajectory moves rapidly towards the $c$-nullcline; it does not ensure that the trajectory stays close to the $c$-nullcline as $T_z\to \infty$. The long-time validity of the zeroth-order approximation is regulated by $\varepsilon_L$, and $\varepsilon_D \leq \varepsilon_L$. Hence, $\varepsilon_L$ is \textit{the} fundamental small parameter, since it is independent of a particular choice of initial data and provides an upper bound on the $\limsup$:
\begin{equation}
    \displaystyle \limsup_{t\to \infty} (c-c_0(p))^2 \leq \lambda^2\varepsilon_L^2.
\end{equation}
For example, $\varepsilon_D$ can be made arbitrarily small by taking $s_0$ to be arbitrarily large, but the accuracy of the approximation in regimes were $s_0$ is large may only be temporary at best (see {\sc figure \ref{fig:final}}). Finally, the magnitude of $\varepsilon_L$ is regulated by the TFPs. The parameter $\varepsilon_L$ determines ``how small" a TFP must be in comparison to other parameters in order to ensure the validity of the tQSSA. Thus, there is a fundamental ``hierarchy" or ``ordering" of parameters:
\begin{equation}\label{eq:ordering}
    \text{TFP} \prec \varepsilon_L \prec \varepsilon_D \prec \varepsilon_T.
\end{equation}
From the perspective of the ordering given in (\ref{eq:ordering}), timescale separation should be seen as an effect, not a cause. Segel and Slemrod~\cite{Segel1989} understood this; they were very clear in issuing their warning of sufficient versus necessary conditions for the validity of QSS reductions, and stressed the fact that timescale separation is merely a necessary condition. 

The takeaway from our analysis is that timescale separation is really a local metric when it comes to assessing the validity of a QSS reduction. This is primarily because timescales are dependent on initial conditions, and the vector field can change dramatically once the trajectory leaves a small neighborhood containing the initial conditions.
Roussel and Fraser~\cite{Roussel:1990:GSSA} seem to have been the first to really understand this in the context of the Michaelis--Menten reaction mechanism (\ref{eq:mm1}), as their work focused on approximating the actual invariant slow manifold. In fact, to the best of the authors' knowledge, Roussel and Fraser (see~\cite{RousselThesis}, pp. 42-47) were the first to provide a solid geometrical argument that $\eta$ is in fact the ``global"  qualifier for the validity of the sQSSA instead of $\varepsilon_{SS}$. By approximating the slow manifold to a high degree of accuracy, Roussel and Fraser~\cite{Fraser1988,Roussel:1991:ASS} were able to demonstrate that $\eta \ll 1$ is necessary for the long-time validity of the sQSSA. 

\subsection{Tightening the limit supremum}
It is worth pointing out that a long-time bound smaller than $\varepsilon_L$ can be obtained when $e_0 < s_0$. This is to be expected, since $\underline{\varepsilon}$ does not vanish if $s_0>e_0$ as $K_M \to 0$. To obtain a tighter bound, note that the upper bound of $c-h^+(p)$ given in (\ref{eq:UP}) is rather liberal. What we want to do is compute $\max (c-h^+(p))$ with more care. First, since $c<h^+(p)$ and $0<h^+(p)$, it follows that $c-h^+(p) <0$. Second, $h^+(p)$ is a decreasing function of $p$, which means $\min h^+(p)$ occurs when $p=s_0$. However, if $c \neq 0$, then $p \leq s_0-c$. Thus, we want to maximize:
\begin{equation}
    c-h^+(s_0-c) \equiv \theta(c) < 0. 
\end{equation}
Since $\theta'(c) >0$ and $\theta(c) <0$, the maximum of $\theta(c)$ occurs when $c=\lambda$, from which it is straightforward to show 
\begin{equation}
 \cfrac{\text{d}\mathcal{E}^2}{\text{d}t}   \leq -\mathcal{E}^2k_1|\theta(\lambda)|+\cfrac{\bigg(\max\bigg|\cfrac{\text{d}h^-}{\text{d}p}\bigg|\bigg)^2\bigg(\sup\bigg|\cfrac{\text{d}p}{\text{d}t}\bigg|\bigg)^2}{k_1|\theta(\lambda)|},
\end{equation}
where $\theta(\lambda)$ is given by:
\begin{equation}
    \theta(\lambda) = \cfrac{\lambda}{2} - \cfrac{1}{2}\bigg(e_0+K_M+\sqrt{(e_0+K_M+\lambda)^2-4e_0\lambda}\bigg).
\end{equation}
An even tighter bound on $\theta(c)$ can be obtained by noting that $\lambda \leq e_0$, and thus $\theta(\lambda) \leq \theta(e_0)$:
\begin{equation}
    \theta(c) \leq \theta(e_0) = -\cfrac{1}{2}\bigg(K_M+\sqrt{4e_0K_M+K_M^2}\bigg).
\end{equation}
The application of Gronwall's Lemma reveals
\begin{equation}
    \limsup_{t\to \infty} (c-h^-(p))^2 \leq \lambda^2 \bigg(\cfrac{\eta}{1+\eta}\bigg)^2\bigg(\cfrac{1}{1+\kappa}\bigg)^2 \bigg(\cfrac{2K_M}{K_M+\sqrt{K_M^2+4e_0K_M}}\bigg)^2  \equiv \lambda^2 \varepsilon_{LT}^2,
\end{equation}
which vanishes as $K_M \to 0$. Thus, the long-time bound ``$\varepsilon_{LT}$" vanishes if $e_0\to 0$, $K_M\to 0$, $e_0\to \infty$, $K_M \to \infty$, or $k_2 \to 0$. Finally, note that
\begin{equation}
    \cfrac{K_M}{e_0-\lambda +K_M} = \cfrac{K_M}{-\theta(e_0)} = \cfrac{K_M}{|\theta(e_0)|}
\end{equation}
when $e_0=s_0$, so there is a nice ``overlap" of small parameters. Note also that the rate of convergence as $K_M\to 0$ is $\mathcal{O}(\sqrt{K_M})$ when $e_0 \leq s_0$. This is as expected. In fast/slow systems that contain a transcritical point, trajectories remain at a distance that is $\mathcal{O}(\sqrt{\varepsilon})$ from the transcritical point once $0<\varepsilon\ll 1$ (again, see~\cite{kuehn2015multiple} for details).

\section{Discussion} \label{sec:discussion}
The primary goal of this work was to illustrate how a combination of 
scaling and non-scaling methods can be employed to determine 
the validity of the sQSSA, extended QSSA, rQSSA, and tQSSA. 
Collectively, our work 
provides an analysis that admits a considerably clearer understanding 
of each QSSA, its corresponding origin, and its validity. Using energy methods, we have recovered a clear set of associated ``small parameters'' that correspond 
to the three types of fast/slow dynamics considered under the
QSSA: the sQSSA, the rQSSA, and the tQSSA  (see 
{\sc Table \ref{table:1}}). \\
\begin{table}[htb!]
\centering
\begin{tabular}{l*{1}|{c}|{c}}
\textbf{Parameter}            & \textbf{Approximation} &  \textbf{Critical Manifold (Set)} \\
\hline
\hline
\;\;\quad $\eta$ & sQSSA & $1-\bar{s}-\bar{p}=0$  \\
\;\;\quad $\underline{\varepsilon}$           &  rQSSA & $\{(\hat{c},\bar{p})\in \mathbb{R}^+|\; (1-\hat{c})(1-\hat{c}-\bar{p})=0\}$ \\
\;\;\quad $\nu$          & tQSSA & $\bar{c}_0(\bar{p};K_S)$ \\
\hline
\end{tabular}
\center
\caption{{\bf Small parameters and critical sets for the standard, 
reverse, and total QSSA}. Each of the critical manifolds that appear 
when either $\eta, \underline{\varepsilon}$ or $\nu$ are identically 
zero correspond to zero product formation. Hence, the tQSSA is valid 
in all three cases. The notation $\bar{c}_0(\bar{p};K_S)$ denotes the $c$-nullcline with $K_M$ replaced with $K_S$, which occurs when $k_2=0.$ }\label{table:1}
\end{table}

The closure of the rQSSA problem for the MM reaction mechanism (\ref{eq:mm1}) is a testament to the power of using a combination of scaling and non-scaling methods, as we were able to derive a ``small parameter'' that ensures 
the validity of the rQSSA, $\underline{\varepsilon} \equiv K_M/(e_0-\lambda)$. In previous studies \cite{Segel1989,SCHNELL2000}, it was assumed that the rQSSA is valid
only when $s_0 \ll e_0$. The small parameter, $\underline{\varepsilon}$, is a major improvement on the small parameters previously reported in the 
literature, because it quantifies and illustrates 
the validity of the rQSSA in regimes where $e_0 \lesssim s_0$. Moreover, for the first time, we have shown that 
the critical manifold associated with the rQSSA at finite initial 
enzyme concentration 
contains a transcritical singularity that significantly influences 
the validity of the rQSSA whenever $e_0 \leq s_0$ and $K_M\ll e_0$, since 
the critical (sub) manifolds undergo an exchange of stability at 
the singular point. On the other hand, this singularity is removed from the physical domain of interest 
if one takes $e_0\to \infty$. This is the distinguishing difference 
between our analysis and the previous analysis of Schnell and 
Maini~\cite{SCHNELL2000}.

Our analysis has practical implications in the laboratory for
the use of mathematical approximations under the rQSSA. Experimental assays are typically implemented to quantify the enzyme activity of a specific reaction. 
This is done by fitting experimental data to a particular model 
equation derived from a QSS reduction (i.e., either the sQSSA, the rQSSA, or the tQSSA). When the rQSSA is valid, progress 
curve experiments for $p$ can be fitted to the product
formation rate expression~(\ref{eq:REA}), which will allow
the estimation of the catalytic rate constant, $k_2$. In mathematical terms, the procedure 
of fitting experimental data to a mathematical model so as to estimate kinetic parameters is an \textit{inverse problem}~\cite{2016-Stroberg-BPC,KIM}. Contrary 
to what has been reported, experimental assays can be carried out 
when initial substrate and enzyme ratios are of similar magnitude,
and the rQSSA can be used to estimate $k_2$, as long as $\underline{\varepsilon}\ll 1.$  

The mathematical analysis of the Michaelis-Menten reaction mechanism (\ref{eq:mm1}) has a long and rich history. We believe that our analysis has deepened the overall understanding
of the QSSA, and that our approach provides an important improvement to the 
implementation of this approximation in enzyme kinetics. The methods we employed can be adapted to analyze 
complex enzyme catalyzed reactions, or multiscale dynamical systems. This will be the focus of our future work. 

\section*{Acknowledgement}
We are grateful to Sebastian Walcher (RWTH Aachen University) and Wylie Stroberg (University of Michigan Medical School) for reading 
and providing comments to previous versions of this manuscript. Justin Eilertsen is supported by the University of Michigan 
Postdoctoral Pediatric Endocrinology and Diabetes Training Program
``Developmental Origins of Metabolic Disorder'' (NIH/NIDDK Grant: 
K12 DK071212).

\appendix 

\section{Appendix} \label{sec:appen}
Here we present a detailed outline of the the method of slowly-varying Lyapunov functions applied to the governing equations of the 
Michaelis--Menten reaction mechanism as discussed in 
Section~\ref{sec:enslavement}

\subsection{Gronwall's lemma}

This differential version can be found in~\cite{multiscale}. Suppose that $u(t)\in C^1([0,T];\mathbb{R}_0^+)$ satisfies\footnote{The notation $L^1(\Omega;\mathbb{R}_0^+)$ refers to the space of functions that are Lebesgue-integrable on $\Omega$ that map $\Omega$ to $\mathbb{R}_0^+$.}
\begin{equation}
    \cfrac{\text{d}u}{\text{d}t} \leq \lambda u + f(t), \qquad u(0)=u_0.
\end{equation}
Then, for $\lambda \in \mathbb{R}$ and $f(t)\in L^1([0,T];\mathbb{R}_0^+)$, it holds that:
\begin{equation}
    u(t) \leq \displaystyle e^{\displaystyle \lambda t}u_0 + e^{\displaystyle \lambda t}\int_0^t e^{\displaystyle -\lambda \tau}f(\tau)\text{d}\tau, \quad \forall t \in[0,T].
\end{equation}

\subsection{The ``\textit{contractivity}" condition: Existence of a Lyapunov function}
It will help to briefly review the concept of a Lyapunov function before we introduce the idea of a ``contractivity" condition. Consider a simple dynamical system of the form
\begin{equation}\label{eq:dyA}
x' = f(x),
\end{equation}
where ``$'$" denotes differentiation with respect to time, $t$. The vector field $f:\mathbb{R}^n\mapsto \mathbb{R}^n$ \textit{generates} a flow, and any scalar-valued function of $x$ (we will call this function $L(x): \mathbb{R}^n\mapsto \mathbb{R}$) will change with respect to time according to:
\begin{equation}
    L' = f(x) \cdot \nabla L \equiv \mathcal{G}L,
\end{equation}
where the operator ``$\mathcal{G}$" is known as a ``generator." Now suppose that the dynamical system (\ref{eq:dyA}) has an equilibrium point at $x=x^*$, so that $f(x^*)=0$. Linearization is typically the preferred method of choice to determine the stability of $x^*$. However, it is well-known that if 
\begin{equation}\label{eq:req}
    0< L, \;\;\&\;\;\mathcal{G}L < 0, \quad \forall x\in\mathcal{B}\backslash \{x^*\},
\end{equation}
where $x^*$ is an interior point of the neighborhood (ball) $\mathcal{B}$, then the equilibrium point $x^*$ is asymptotically stable, and ``$L$" is called a \textit{Lyapunov} function.

So what does this have to do with a contraction condition? Keeping the idea of a Lyapunov function in mind, consider a fast/slow system:
\begin{subequations}
\begin{align}
    x' &= \varepsilon f(x,y),\quad x(0)=x_0\label{eq:fA1}\\
    y' &= g(x,y),\quad y(0)=y_0\label{eq:fA2},
\end{align}
\end{subequations}
where the slow variable, $x\in\mathbb{R}^n$, the fast variable $y\in\mathbb{R}^m$, and thus $f(x,y):\mathbb{R}^n\times \mathbb{R}^m \mapsto \mathbb{R}^n$ and $g(x,y):\mathbb{R}^n\times \mathbb{R}^m \mapsto \mathbb{R}^m$. The ``fast subsystem" associated with (\ref{eq:fA1})--(\ref{eq:fA2}) is
\begin{equation}
y' = g(x_0,y).
\end{equation}
Assume the fast subsystem has an equilibrium point $y^* = h_0(x_0)$, where $h:\mathbb{R}^n \mapsto \mathbb{R}^m$, and $g(x_0,h_0(x_0))=0$. Since $x_0$ is a constant with respect to the fast subsystem, $g$ can be treated as a map from $\mathbb{R}^m$ to $\mathbb{R}^m$. Once again, the stability of $y^*$ can be found via linearization. However, suppose we want to determine an appropriate Lyapunov function. Generally, finding a Lyapunov function can be difficult, and so it is common practice to write down a simple polynomial function and ``check" to see if it satisfies the requirements given in $(\ref{eq:req})$. The simplest choice is a quadratic function, $L(y) =||y-h_0(x_0)||^2$, which is locally positive-definite since $L(y)>0, \;\; \forall y\neq y^*=h_0(x_0)$. The second property we  need to check is whether or not $L'$ is locally negative-definite; thus, we need to show that
\begin{equation}\label{eq:ND}
L' = 2\langle(y-y^*),g(x_0,y)\rangle < 0,\quad \forall y\in \mathcal{B}\backslash \{y^*\},
\end{equation}
where $\langle,\rangle$ denotes the usual scalar product between two vectors, and $||v|| =\sqrt{\langle v,v \rangle}$. It is not entirely clear how go about showing that $L'$ is negative definite, since $g(x_0,y)$ and $(y-y^*)$ can change signs. However, a little ``trick" can go a long way. First, rewrite (\ref{eq:ND}) as
\begin{equation}
    2\langle (y-y^*),g(x_0,y)\rangle = 2\langle (y-y^*), g(x_0,y)-g(x_0,y^*)\rangle.
\end{equation}
Now suppose we know that there is some positive number ``$\zeta$" such that
\begin{equation}\label{eq:Acon}
    \langle(y-y^*),g(x_0,y)-g(x_0,y^*)\rangle \leq -\zeta||y-y^*||^2.
\end{equation}
If such a number ``$\zeta$" exists, then we can bound $L'$
\begin{equation}\label{eq:diffINQ}
    L' = 2\langle(y-y^*),g(x_0,y)\rangle \leq -2\zeta ||y-y^*||^2 \leq 0,
\end{equation}
and thus we have shown that $L'$ is locally negative-definite. Consequently, $L=||y-y^*||^2 \equiv ||y-h_0(x_0)||^2$ is a Lyapunov function. Why then, do we refer to  (\ref{eq:Acon}) as a ``contractivity" condition? Well, notice that (\ref{eq:diffINQ}) defines a differential inequality
\begin{equation}
L' \leq -2\zeta L,
\end{equation}
and integrating this inequality (courtesy of Gronwall's lemma) yields
\begin{equation}
    L\leq L(0)e^{\displaystyle - 2\zeta t}.
\end{equation}
Consequently, the distance between $y_0$ and $h_0(x_0)$ decays to zero exponentially; hence, the distance ``contracts."

In the subsection that follows, we will investigate the growth of $L = ||y-h_0(x)||^2$, but we will allow $x$ to change in time, since $x$ will vary slowly in comparison to $y$ when $0< \varepsilon \ll 1$. Thus, we will be interested in finding an upper bound on 
\begin{equation}
\cfrac{\text{d}L}{\text{d}T} \equiv \dot{L},
\end{equation}
where $T =\varepsilon t$ is the slow timescale. The motivation for doing this should be obvious: we want to compute an estimate on how well the reduced system, $\dot{x} = f(x,h_0(x))$, approximates the full system. Because we will take into account the temporal variation of $x$, $L$ is no longer a Lyapunov function. However, because $x$ varies slowly, we have chosen to refer to $L$ in this case as a \textit{slowly-varying Lyapunov function}. One could also refer to this procedure as simply an ``energy method," since we are ultimately bounding the energy (norm) of the error $y-h_0(x)$.
\subsection{General procedure: Estimation of bounds via energy methods}
Let us start by estimating bounds on the enslavement of the 
fast variable. The calculation outlined below follows from the original work of Tikhonov \cite{tikhonov1952} and Gradshtein \cite{Gradshtein1953}, and the particular calculation we outline here can be found in \cite{multiscale}. This procedure will then be used to develop bounds 
for enslavement in the $(p,c)$ coordinate system. We will express the procedure in terms with respect to a the fast/slow system (\ref{eq:fA1})--(\ref{eq:fA2}) and, since our specific problem of interest is two-dimensional, we will assume $x\in\mathbb{R}$ and $y\in\mathbb{R}$. However, the reader should keep in mind that this procedure extends to higher-dimensional problems. 

With respect to the slow time, the system (\ref{eq:fA1})--(\ref{eq:fA2}) is expressed in form
\begin{subequations}
\begin{align}
\dot{x} &= f(x,y),\\
\varepsilon \dot{y} &= g(x,y),
\end{align}
\end{subequations}
where ``$\dot{\phantom{x}}$" denotes differentiation with respect to the slow time, $T$. Let $g(\cdot,y)$ be contractive in the sense that
\begin{equation}\label{eq:contractivity}
 (g(x,y_1)-g(x,y_2))\cdot(y_1-y_2) \leq -\zeta \cdot (y_1-y_2)^2,
\end{equation}
where the ``$\cdot$" in (\ref{eq:contractivity}) and in what follows denotes scalar multiplication. Next, suppose that $h_0(x)$ satisfies $g(x,h_0(x))=0$. Defining 
$y=h_0(x) +z$, where $z=y-h_0(x)$, it follows that
\begin{equation*}
\dot{z}  = \dot{y} - D_x h_0(x) \cdot \dot{x}. 
\end{equation*}
Then, the differential equation for the energy, $L\equiv z ^2=(y-h_0(x))^2$, is
\begin{equation}
\cfrac{1}{2}\dot{z^2} = \cfrac{1}{\varepsilon} z\cdot g(x,z+h_0(x)) -z\cdot D_x h_0(x)\cdot \dot{x}.
\end{equation}
Next, we want to get this expression into a useful form so that we can apply Gronwall's lemma. There are essentially two ``tricks" that can be utilized to manipulate this expression: (i.) we can ``add/subtract" zero, or (ii.) we can ``multiply" by one. We will choose (i). The term $\varepsilon^{-1} z \cdot g(x,z+h_0(x))$ can be expanded by ``subtracting" zero,
\begin{equation}
    \cfrac{1}{\varepsilon}z\cdot g(x,z+h_0(x)) = \cfrac{1}{\varepsilon}\bigg[ z\cdot g(x,z+h_0(x)) -  z \cdot g(x,h_0(x))\bigg],
\end{equation}
which holds since $g(x,h_0(x))=0$. Thus, we obtain
\begin{equation}
\cfrac{1}{2}\dot{z^2}= \cfrac{1}{\varepsilon}\bigg[ z \cdot g(x,z+h_0(x))- z \cdot g(x,h_0(x))\bigg] -z \cdot D_x h_0(x) \cdot \dot{x}.\label{eq:Lbound}
\end{equation}
Our primary interest is not in finding an actual solution to (\ref{eq:Lbound}). Instead, we want to compute an upper bound on the growth of $z^2$, which means we need to convert the differential equation (\ref{eq:Lbound}) into a differential \textit{inequality}. To do this, we will take advantage of the ``contractivity" condition. First rewrite the first term on the right hand side of (\ref{eq:Lbound}) as:
\begin{equation}
    \bigg[z \cdot g(x,z+h_0(x)) -  z \cdot g(x,h_0(x))\bigg] =  \bigg[g(x,z+h_0(x))-g(x,h_0(x))\bigg]\cdot \bigg[z+h_0(x)-h_0(x)\bigg].
\end{equation}
Let $z+h_0(x)$ denote ``$y_1$" in (\ref{eq:contractivity}), and let $h_0(x)$ denote ``$y_2$" in (\ref{eq:contractivity}). Then,
\begin{equation}\label{eq:cont}
(g(x,z+h_0(x))-g(x,h_0(x)))\cdot (z+h_0(x)-h_0(x)) \leq -\zeta \cdot z^2
\end{equation}
It follows from (\ref{eq:cont}) that
\begin{equation}\label{eq:diffEQ}
\cfrac{1}{2}\dot{z^2}\leq -\cfrac{\zeta}{\varepsilon} z^2 - z \cdot D_x h_0(x) \cdot \dot{x}.
\end{equation}
Maximizing the inequality (\ref{eq:diffEQ}) yields
\begin{equation}\label{eq:ding}
\cfrac{1}{2}\dot{z^2} \leq -\cfrac{\zeta}{\varepsilon} z^2 + |z|\max|D_x h_0(x)|\max|\dot{x}|.
\end{equation}
In order to use Gronwall's lemma, we want the right hand side of (\ref{eq:ding}) to be in terms of $z^2$ instead of $|z|$. To do this, 
let $a \equiv |z|$ and $b \equiv \max|D_x h_0(x)|\max|\dot{x}|$. 
Then, from Cauchy's inequality, ($a\cdot b \leq \delta a^2 + b^2/4\delta $), 
it holds that
\begin{equation}
|z|\max|D_x h_0(x)|\max|\dot{x}| \leq \delta z^2 + \cfrac{(\max|D_x h_0(x)|\max|\dot{x}|)^2}{4\delta}, \quad \forall \delta>0.
\end{equation}
Choosing $\delta= \zeta/2\varepsilon $ yields
\begin{equation}\label{eq:inter}
|z|\max|D_x h_0(x)|\max|\dot{x}| \leq \cfrac{\zeta}{2\varepsilon}z^2 +\varepsilon \cfrac{(\max|D_x h_0(x)|\max|\dot{x}|)^2}{2\zeta}.
\end{equation}
Substitution of (\ref{eq:inter}) into (\ref{eq:ding}) reveals the differential inequality
\begin{equation}\label{eq:gronwall}
\dot{z^2} \leq -\cfrac{\zeta}{\varepsilon}z^2 +\varepsilon \cfrac{(\max|D_x h_0(x)|\max|\dot{x}|)^2}{2\zeta}.
\end{equation}
The expression (\ref{eq:gronwall}) is now in a form that be integrated directly, and it follows from Gronwall's lemma that,
\begin{subequations}
\begin{align}
z^2 &\leq z^2(0)\displaystyle e^{\displaystyle -\zeta T/\varepsilon} + \varepsilon^2 \cfrac{(\max|D_x h_0(x)|\max|\dot{x}|)^2}{\zeta^2}\bigg(1-\displaystyle e^{\displaystyle -\zeta T/\varepsilon}\bigg),\\
&\leq z^2(0)\displaystyle e^{\displaystyle -\zeta T/\varepsilon} + \varepsilon^2 \cfrac{(\max|D_x h_0(x)|\max|\dot{x}|)^2}{\zeta^2}.
\end{align}
\end{subequations}
Finally, from the triangle inequality, we obtain:
\begin{equation}\label{eq:bounder}
|z| \equiv |y-h_0(x)| \leq |y_0-h_0(x_0)|\displaystyle e^{\displaystyle -\zeta T/2\varepsilon} + \varepsilon \cfrac{\max|D_x h_0(x)|\max|\dot{x}|}{\zeta}.
\end{equation}



\subsection{Computation of the upper bound for enslavement of \texorpdfstring{$c$}{Lg} in \texorpdfstring{$(p,c)$}{Lg} coordinates}
Here we will directly calculate the constants $\zeta, \max |D_x h(x)|$ and $\max |\dot{x}|$ need for the upper bound. Note that the computation of $\zeta$ is calculated over the domain $(p,c)\in \Omega \equiv [0,s_0] \times [0,\lambda]$, since trajectories that start on $(p,c)(0)=(0,0)$ are confined to this domain due to conservation. 

In the $(p,c)$ coordinate system, the $c$-nullcline is given by:
\begin{equation}
h(p) = \cfrac{1}{2}\bigg(e_0+K_M+s_0-p - \sqrt{(e_0+K_M+s_0-p)^2-4e_0(s_0-p)}\bigg).
\end{equation}
Finding $\max|h'(p)|$ is relatively straightforward. First,
\begin{equation}
\cfrac{\text{d}^2h}{\text{d}p^2} =-\cfrac{2e_0K_M}{((e_0+K_M+s_0-p)^2-4e_0(s_0-p))^{3/2}} < 0,
\end{equation}
and therefore the derivative is monotonically decreasing on the interval $p\in[0,s_0]$. Thus, we obtain:
\begin{equation}
    \max \bigg| \cfrac{\text{d}h}{\text{d}p}\bigg| = -\cfrac{\text{d}h}{\text{d}p}\bigg|_{p=s_0}= \cfrac{e_0}{K_M+e_0}.
\end{equation}
For contraction, we define
\begin{equation}
G(c,p) \equiv \dot{c} = k_1c^2 -k_1(K_M+e_0+s_0-p)c + k_1e_0(s_0-p).
\end{equation}
We will use the Mean Value Theorem to find $\zeta$. The derivative of $G(c,p)$ is given by
\begin{equation}
\cfrac{\partial G}{\partial c} = -k_1(e_0+s_0+K_M-2c-p) \leq -k_1(e_0+K_M+s_0) + \sup k_1(2c+p).
\end{equation}
We need to find $\sup (2c+p)$. A straightforward bound is $2\lambda +s_0$. However, we want to find the supremum subject to the constraints: $s_0 = s+c+p$ and $c \leq \lambda$. Employing the conservation law ``$s_0=s+c+p$ yields
\begin{equation}
    2c+p = s_0-s+c \leq s_0 + \lambda.
\end{equation}
Consequently, since $2c +p = s_0+c-s \leq s_0 +\lambda$, it follows that
\begin{equation}\label{eq:EQQ}
|c-c_0(p)| \leq \lambda e^{\displaystyle-t\zeta_T/2} + \cfrac{\eta\lambda}{(1+\eta)(1+\kappa)}\bigg(\cfrac{K_M}{e_0-\lambda + K_M}\bigg),
\end{equation}
where $\zeta_T$ is given by
\begin{equation}
\zeta_T = k_1(e_0-\lambda+K_M).
\end{equation}
Diving both sides of (\ref{eq:EQQ}) by $\lambda$ and expressing $p=s_0\bar{p}$ yields the dimensionless form of the inequality
\begin{equation}
|\underline{c}-\bar{c_0}(\bar{p})| \leq e^{\displaystyle-\tilde{T}/2\underline{\varepsilon}\nu} + \cfrac{\underline{\varepsilon}}{(1+\underline{\varepsilon})}\cfrac{\eta}{(1+\eta)}\cfrac{1}{(1+\kappa)}.
\end{equation}
Since $\nu \leq 1$, we obtain
\begin{equation}
|\underline{c}-\bar{c_0}(\bar{p})| \leq e^{\displaystyle-\tilde{T}/2\underline{\varepsilon}} + \underline{\varepsilon}.
\end{equation}
Furthermore, since $\zeta_T \geq k_1K_M$, it holds that
\begin{equation}
|\underline{c}-\bar{c_0}(\bar{p})| \leq e^{\displaystyle-\bar{T}/2\varepsilon_{SS}} + \eta.
\end{equation}
Let us define a ``long-time epsilon, $\varepsilon_{L}$,"
\begin{equation}
    \varepsilon_{L} \equiv \cfrac{\underline{\varepsilon}}{(1+\underline{\varepsilon})}\cfrac{\eta}{(1+\eta)}\cfrac{1}{(1+\kappa)}.
\end{equation}
Expressing the exponential in terms of Tzafriri's timescale, $T_z = k_2\lambda t/s_0$, we obtain:
\begin{equation}
|\underline{c}-\bar{c_0}(\bar{p})| \leq e^{\displaystyle-s_0T_z/2\lambda \underline{\varepsilon}\nu} + \varepsilon_{L}.
\end{equation}
Moreover, we can define a ``decay epsilon," $\varepsilon_{D}\equiv s_0^{-1}\lambda\underline{\varepsilon}
\nu$, and write
\begin{equation}
|\underline{c}-\bar{c_0}(\bar{p})| \leq e^{\displaystyle-T_z/2\varepsilon_D} + \varepsilon_{L},
\end{equation}
where $\varepsilon_{D}\leq \varepsilon_L.$

\subsection{The tQSSA as used in practice}

The $c$-nullcline is rarely used for parameter estimation in progress curve experiments. Typically, at high enzyme concentrations, one uses
\begin{equation}
c\approx \cfrac{e_0(s_0-p)}{e_0+K_M+s_0-p}
\end{equation}
in order to employ
\begin{equation}
    \cfrac{\text{d}p}{\text{d}t} \approx \cfrac{k_2e_0(s_0-p)}{e_0+K_M+s_0-p}
\end{equation}
for parameter estimation (see \cite{KIM} for a good discussion on the tQSSA and parameter estimation). Calculating an upper bound on 
\begin{equation}
    |\mathcal{E}_T| \equiv \bigg| c- \cfrac{e_0(s_0-p)}{e_0+K_M+s_0-p} \bigg|
\end{equation}
is done using the methods previously. First,
\begin{subequations}
\begin{align}
    \cfrac{1}{2}\cfrac{\text{d}\mathcal{E}_T^2}{\text{d}t} &= -k_1(e_0+K_M+s_0-p)\mathcal{E}_T^2 + \bigg[k_1c^2 + \cfrac{k_2c e_0(e_0+K_M)}{(e_0+K_M+s_0-p)^2}\bigg]\mathcal{E}_T,\\
    &\leq -k_1(e_0+K_M)\mathcal{E}_T^2 + \bigg[k_1\lambda^2 + \cfrac{k_2\lambda e_0}{e_0+K_M}\bigg]|\mathcal{E}_T|.
    \end{align}
\end{subequations}
Applying the Cauchy-$\delta$ inequality with $\delta = k_1(e_0+K_M)/2$, followed by integration, yields
\begin{equation}
\mathcal{E}_T^2 \leq \mathcal{E}_T^2(0) \displaystyle e^{\displaystyle -k_1(K_M+e_0)t} + \cfrac{1}{(k_1(e_0+K_M))^2}\bigg[k_1\lambda^2 + \cfrac{k_2\lambda e_0}{e_0+K_M}\bigg]^2.
\end{equation}
Finally, after scaling and algebraic simplification, it holds that:
\begin{equation}\label{eq:pTQ}
  |\mathcal{E}_T| \leq |\mathcal{E}_T(0)|e^{\displaystyle -k_1(K_M+e_0)t/2} + \lambda \bigg[\cfrac{\lambda}{e_0+K_M} + \nu\bigg(\cfrac{ e_0}{e_0+K_M}\bigg)\bigg(\cfrac{K_M}{e_0+K_M}\bigg)\bigg].
\end{equation}
In terms of $T_z$, the inequality (\ref{eq:pTQ}) is
\begin{subequations}
\begin{align}
   |\bar{\mathcal{E}}_T| &\leq |\bar{\mathcal{E}}_T(0)|e^{\displaystyle -s_0T_z/2\lambda \tilde{\epsilon}} +  \cfrac{\eta(1+\nu\tilde{\epsilon})}{1+\eta},\qquad \qquad \;\;\text{if}\;\;e_0 \leq s_0,\\
   |\bar{\mathcal{E}}_T| &\leq |\bar{\mathcal{E}}_T(0)|e^{\displaystyle -s_0T_z/2\lambda \tilde{\epsilon}} +  \cfrac{\lambda}{e_0+K_M}+ \eta\nu\tilde{\epsilon},\qquad \text{if}\;\;s_0 < e_0
   \end{align}
\end{subequations}
where $\bar{\mathcal{E}}_T =\lambda^{-1}\mathcal{E}_T$ and $\tilde{\epsilon}\equiv K_M/(e_0+K_M)$. Thus, the tQSSA as used in practice is valid at very low $e_0$ as well as very high $e_0$.












\end{document}